\documentclass[11pt,a4paper]{article}
\usepackage[english]{babel}
\usepackage[utf8]{inputenc}
\usepackage[T1]{fontenc}
\usepackage{lmodern,amsmath,amssymb,amsthm,mathtools,mathrsfs}
\usepackage[a4paper,margin=2.35cm]{geometry}
\usepackage{enumitem,hyperref,xcolor,tikz}
\hypersetup{colorlinks=true,linkcolor=blue!50!black,urlcolor=blue!50!black}
\newtheorem{theorem}{Theorem}[section]
\newtheorem{lemma}[theorem]{Lemma}
\newtheorem{proposition}[theorem]{Proposition}
\newtheorem{corollary}[theorem]{Corollary}
\theoremstyle{definition}\newtheorem{definition}[theorem]{Definition}
\theoremstyle{remark}\newtheorem{remark}[theorem]{Remark}
\newcommand{\N}{\mathbb N}
\newcommand{\Lip}{\operatorname{Lip}}
\newcommand{\diam}{\operatorname{diam}}
\newcommand{\dist}{\operatorname{dist}}

\newcommand{\JL}{\mathscr{JL}}
\title{Separable $\mathcal L_{\infty}$-spaces with Kalton's property $Q$}
\author{Jesús Suárez de la Fuente}
\begin{document}
\maketitle
\begingroup
\renewcommand{\thefootnote}{}
\endgroup
\begin{abstract}
We construct separable $\mathcal L_\infty$-spaces with Kalton's property $Q$,
which, to the best of our knowledge, are the first known examples of this kind.
Our technique employs the Bourgain--Pisier construction and a quantitative
analysis of metric types in limits of finite-dimensional
$\eta$-admissible embeddings. The argument rests on two complementary
estimates: approximate radial centres yield asymptotic separation of tree
increments, while an admissibility inequality yields asymptotic gain,
leading to concentration on
Kalton's interlacing graphs. This leads to an abstract criterion
for property $Q$ in terms of asymptotic separation $AS_\theta$
and asymptotic gain $AG_\alpha$. We illustrate the limitations of this
criterion with Johnson--Lindenstrauss spaces. Our examples exclude uniform embeddings of the unit ball of $c_0$,
strengthening in a nonlinear direction the classical Bourgain--Delbaen
phenomenon. As an application, every separable Banach space embeds
isometrically into a separable $\mathcal L_\infty$-space whose unit ball
is not a uniform retract of the unit ball of its bidual.
\end{abstract}
\setcounter{tocdepth}{2}
\tableofcontents
\clearpage
\section{Introduction}

Kalton's property $Q$ is a concentration principle for Lipschitz maps
on a family of combinatorial graphs. Its significance lies in a
uniformity requirement: the concentration estimate must be independent
of the length of the tuples indexing the graph. This makes property
$Q$ sensitive to the large-scale metric geometry of a Banach space,
even when the maps under consideration have bounded range.

To recall the definition, for $k\ge1$ let $[\mathbb N]^k$ denote the
set of strictly increasing $k$-tuples
$\bar n=(n_1,\ldots,n_k)$ of positive integers. Two distinct tuples
$\bar n,\bar m\in[\mathbb N]^k$ are joined by an edge when
\begin{equation}\label{eq:intro-interlacing}
 n_1\le m_1\le n_2\le m_2\le\cdots\le n_k\le m_k,
 \qquad\text{or the same inequalities with $\bar n$ and $\bar m$ exchanged}.
\end{equation}
The corresponding \emph{interlacing graph} is equipped with its
shortest-path distance $d_I$. If $M\subset\mathbb N$ is infinite,
$[M]^k$ denotes the tuples with all their entries in $M$.
A metric space $(Z,d)$ has \emph{Kalton's property $Q$} if there exists
$C<\infty$ such that, for every $k\ge1$ and every bounded Lipschitz
map $f:([\mathbb N]^k,d_I)\to Z$, some infinite $M\subset\mathbb N$
satisfies
\begin{equation}\label{eq:intro-Q}
 \operatorname{diam} f([M]^k)\le C\operatorname{Lip}(f).
\end{equation}
We write $Q_Z$ for the infimum of all such constants $C$. In
particular, the same $C$ must work for \emph{every} $k$ and every
$f$, although the set $M$ may depend on both. This formulation
of property $Q$, and its relation to Kalton's original one, is
recalled in the preliminaries; see also \cite{K,Qstab}.

Property $Q$ was introduced by Kalton in his study of nonlinear
embeddings of Banach spaces \cite{K}. His concentration theorem \cite[Theorem~4.1]{K} is
stronger than the property-$Q$ estimate above: in a reflexive space,
the restricted image lies in a ball of controlled radius, not only
in a set of controlled diameter.

\begin{theorem}[Kalton]\label{thm:intro-kalton}
Let $X$ be a reflexive Banach space. For every $k\ge 1$, every
bounded map $f:([\mathbb N]^k,d_I)\to X$, and every $\varepsilon>0$,
there are an infinite set $M\subset\mathbb N$ and a point $x\in X$
such that
\begin{equation}\label{eq:intro-kalton}
 \sup_{\bar n\in[M]^k}\|f(\bar n)-x\|
 \le \operatorname{Lip}(f)+\varepsilon.
\end{equation}
\end{theorem}
Consequently, every reflexive Banach space has property $Q$ and,
with the normalization \eqref{eq:intro-Q}, $Q_X\le 2$.

In sharp contrast, $c_0$ does not have property $Q$.
Indeed, writing $s_n=\sum_{i=1}^n e_i$ for its summing basis, the maps
\[
 f_k(n_1,\ldots,n_k)=\sum_{i=1}^k s_{n_i}
 \qquad (n_1<\cdots<n_k)
\]
are $1$-Lipschitz on the interlacing graphs, while the diameter
of $f_k([M]^k)$ is at least $k$ for every infinite $M\subset\mathbb N$
\cite[Remark~2.7]{Qstab}. Thus the contrast between reflexive
spaces and $c_0$ is already visible on these graphs.

It is natural to ask how this concentration phenomenon interacts
with classes of Banach spaces whose defining features are linear
and finite-dimensional. We address this question for
$\mathcal L_\infty$-spaces: spaces whose finite-dimensional subspaces
can be accommodated, with controlled distortion, in
finite-dimensional spaces resembling $\ell_\infty^n$.
Since $c_0$ is an $\mathcal L_\infty$-space, such local structure
alone does not imply property $Q$.

The main result of this paper gives a contrasting phenomenon:
\emph{there exist separable, infinite-dimensional $\mathcal L_\infty$-spaces
with Kalton's property $Q$.} To the best of our knowledge,
these are the first examples of $\mathcal L_\infty$-spaces known
to have this property. More generally, we obtain these spaces
from Bourgain--Pisier constructions over separable reflexive Banach
spaces. In particular, their unit balls do not admit a uniform
embedding of the unit ball of $c_0$. The point is not merely the
existence of another $\mathcal L_\infty$-space, but the identification
of a geometric mechanism that enforces concentration in a class where
it cannot be inferred from the local linear structure.

This may be viewed as a nonlinear strengthening of a classical
phenomenon. Lindenstrauss and Pe\l{}czy\'nski, and also Lindenstrauss
and Rosenthal, asked whether every infinite-dimensional
$\mathcal L_\infty$-space must contain a subspace isomorphic to $c_0$
\cite{LP1968,LR1969}.
Bourgain and Delbaen answered this question negatively by constructing
separable $\mathcal L_\infty$-spaces without a linear copy of $c_0$
\cite{BD1980}. Our examples go further in a metric direction:
they exclude even a uniform embedding of the unit ball of $c_0$.
The construction employed here belongs to the subsequent
Bourgain--Pisier framework \cite{BP1983,LA}; no assertion is made
that the original Bourgain--Delbaen examples themselves have property $Q$.

Our proof begins with a class of spaces having an admissible
finite-dimensional structure. For $0<\eta<1$, let $S_\eta$ denote
the completion of an increasing union of finite-dimensional spaces
$E_N$ whose inclusions are $\eta$-admissible isometries in the sense
of L\'opez-Abad \cite{LA} (see Section~\ref{subsec:lopez-admissibility} below for the definition).
This class arises naturally in the
Bourgain--Pisier construction. For a separable reflexive space $X$,
the auxiliary construction of L\'opez-Abad provides a Banach space
$\mathcal E$ containing $X\subset\mathcal L_\infty(X)\subset\mathcal E$, with
$\mathcal E/X$ isometric to a space $S_\eta$ of the preceding kind. Hence
$\mathcal L_\infty(X)/X$ embeds isometrically into $S_\eta$.
The admissible finite-dimensional stages are those of $\mathcal E/X$, not
necessarily the stages generating $\mathcal L_\infty(X)$ itself.
We first establish property $Q$ for $S_\eta$, then for
$\mathcal L_\infty(X)/X$, and finally for $\mathcal L_\infty(X)$.

The geometric input is a quantitative consequence of
$\eta$-admissibility. In the push-out description, an admissible
extension of a space $F$ is obtained from $B\oplus_1F$ by
identifying $(t,0)$ with $(0,u(t))$, where $T\subset B$ and
$u:T\to F$ has $\|u\|\le\eta$. Such embeddings retain a quantitative
remnant of the $\ell_1$ geometry of the construction. In particular,
for vectors $z_1,z_2\in S_\eta$ with $z_1+z_2\in E_N$,
\begin{equation}\label{eq:intro-admiss}\tag{A}
 \|z_1\|+\|z_2\|\ge
 \|z_1+z_2\|+(1-\eta)
 \bigl(\operatorname{dist}(z_1,E_N)+
       \operatorname{dist}(z_2,E_N)\bigr).
\end{equation}
This is the admissibility inequality used throughout the proof.
It yields the gain estimate
\begin{equation}\label{eq:intro-acc}
 \|a+w\|\ge\|a\|+2(1-\eta)\operatorname{dist}(w,E_N)
       -\|w\|-2\operatorname{dist}(a,E_N),
\end{equation}
which quantifies the contribution of an increment transverse to a
previous stage.

We next describe the geometric mechanism behind concentration.
Let $S_\eta$ be an abstract limit of finite-dimensional
$\eta$-admissible stages. For a bounded
sequence $(x_j)$ in $S_\eta$ and a free ultrafilter $\mathcal U$,
consider its metric type
\[
 \Phi(z)=\lim_{j,\mathcal U}\|z-x_j\|.
\]
We show that $\Phi$ has an \emph{approximate radial centre}:
there exist $v\in S_\eta$ and $R\ge0$ such that, for every
$z\in S_\eta$,
\begin{equation}\label{eq:intro-centre}
 \|z-v\|+(1-2\eta)R\le\Phi(z)\le\|z-v\|+R.
\end{equation}
The two-sided estimate holds at every point of $S_\eta$; it
is stronger than approximate minimization. Its lower bound
leads to the separation factor
\[
 \theta_\eta=1-2\eta.
\]
Applying these centres recursively to the successors of a finite
tree, and choosing suitable stages along an infinite subset of
$\mathbb N$, we obtain increments that are asymptotically separated
from the stages containing the preceding ones.

\begin{samepage}
More precisely, if $v_s$ denotes the centre associated with a node $s$,
with $v_{\bar n}=f(\bar n)$ at the terminal nodes, and
$u_s=v_s-v_{s^-}$, then
\[
 \boxed{f(\bar n)=v_\varnothing+\sum_{i=1}^k u_{\bar n|i}.}
\]
Thus, the root term disappears when two branches are compared.
\end{samepage}

The second ingredient is asymptotic gain. When two branches of the tree
are strictly interlaced, their increments can be arranged in
alternating order. Applying \eqref{eq:intro-acc} successively shows
that their norms accumulate with a coefficient
\[
 \alpha_\eta=2(1-\eta).
\]
The coefficient in the two-branch estimate is the product of
the asymptotic gain and separation coefficients, minus one:
\begin{equation}\label{eq:intro-factor}
 \lambda_\eta
   =\underbrace{2(1-\eta)}_{\text{gain}}
    \underbrace{(1-2\eta)}_{\text{separation}}-1
   =1-6\eta+4\eta^2.
\end{equation}
For $0<\eta<(3-\sqrt5)/4$, this number is positive. Up to an error
which can be made arbitrarily small after the tuple length $k$ has
been fixed, we obtain for two interlaced leaves
$\bar n,\bar m$ an estimate of the form
\[
 \|f(\bar n)-f(\bar m)\|\ge
 \lambda_\eta\left(
  \sum_{i=1}^k\|u_{\bar n|i}\|+
  \sum_{i=1}^k\|u_{\bar m|i}\|
 \right)-\text{error}.
\]
An application of infinite Ramsey theory makes the norms of the
increments essentially constant at each level of the tree. Since
strictly interlaced leaves are adjacent in the interlacing graph,
the left-hand side is at most $\operatorname{Lip}(f)$. A final
comparison of arbitrary branches gives
\[
 Q_{S_\eta}\le\lambda_\eta^{-1}.
\]
Passing to $\mathcal L_\infty(X)/X$ and then using the
three-space principle for property $Q$ (with the reflexivity of
$X$) yields the announced $\mathcal L_\infty$ examples.

\medskip
\noindent\emph{An abstract separation--gain principle.}
The preceding proof suggests separating the two geometric inputs.
We formulate \emph{asymptotic separation} $AS_\theta$ and
\emph{asymptotic gain} $AG_\alpha$ relative to a common increasing
dense chain of finite-dimensional subspaces. They capture,
respectively, the ability to select transverse increments and the
norm growth that results when these increments are added. Our
abstract criterion states that
\[
 AS_\theta\quad\text{and}\quad AG_\alpha,
 \qquad \alpha\theta>1
 \quad\Longrightarrow\quad
 Q_X\le(\alpha\theta-1)^{-1}.
\]
For the admissible spaces, the two numerical inputs may be taken to be
$\theta=\theta_\eta$ and $\alpha=\alpha_\eta$, which explains
\eqref{eq:intro-factor}. The criterion is sufficient rather than a
characterization. At the end of the paper we examine
Johnson--Lindenstrauss spaces as a test case: although they enjoy
property $Q$, they do not satisfy the stronger asymptotic separation and gain requirements used in this argument. These examples
also clarify the significance of the threshold $1/2$ for the
separation parameter.

A further application concerns uniform retractions of biduals.
Combining the $\mathcal L_\infty$ examples above with Kalton's
renorming construction, we show that every separable Banach space
embeds isometrically into a separable $\mathcal L_\infty$-space whose
unit ball is not a uniformly continuous retract of the unit ball of
its bidual. This addresses a question of Benyamini and Lindenstrauss
\cite[p.~183]{BenyaminiLindenstrauss2000}, also recalled by
Kalton~\cite{Kalton2011}:
\begin{quote}
\emph{Must every Banach space be a Lipschitz retract of its bidual?}
\end{quote}
Kalton constructed a
\emph{nonseparable} Lindenstrauss space, thus an $\mathcal L_\infty$-space,
whose unit ball is not a
uniform retract of its bidual unit ball \cite{Kalton2011}.
More recently, Acuaviva established the failure of such a
retraction for a \emph{separable} space with a Schauder basis
\cite{Acuaviva2026}. Our result provides separable
$\mathcal L_\infty$ counterexamples and, moreover, permits
an isometric copy of any prescribed separable Banach space.

The paper is organized as follows. Section~\ref{sec:prelim}
introduces metric types and interlacing graphs.
Section~\ref{sec:quotient} recalls the admissible Bourgain--Pisier
construction. Sections~\ref{sec:centres}--\ref{sec:two-branches}
develop approximate radial centres, tree extraction and the
two-branch gain estimate. Section~\ref{sec:Q-spaces} proves
property $Q$ for the resulting spaces, while
Section~\ref{sec:kalton-envelope} treats uniform retractions.
Section~\ref{sec:abstract-AS-AG} formulates the abstract
$AS_\theta$--$AG_\alpha$ criterion, and
Section~\ref{sec:JL-AS-AC} examines the Johnson--Lindenstrauss examples.

\section{Preliminaries: metric types and property $Q$}
\label{sec:prelim}

\noindent\textbf{Notation and conventions.}
We use standard terminology from Banach space theory, referring to
\cite{AlbiacKalton,LindenstraussTzafririI} for background. Unless
otherwise stated, $B_X$ denotes the closed unit ball of $X$, and
$\dist(x,E)=\inf_{e\in E}\|x-e\|$. In the Bourgain--Pisier
construction, $X$ will always be separable and infinite-dimensional;
our general statements have precisely the hypotheses specified in
their formulations. We write $\mathcal L_\infty(X)$ for the associated
space and $\mathcal L_{\infty,\lambda}(X)$ when the constant matters.
For property $Q$, metric types and interlacing graphs see
\cite{K,Qstab}; for $\eta$-admissible embeddings and the
Bourgain--Pisier construction see \cite{BP1983,LA}. Throughout the
proofs, $\ell$ denotes a Lipschitz constant and $L$ an infinite
subset of $\mathbb N$.

\subsection{Metric types}

\paragraph{Ultralimits.}
Let $\mathcal U$ be a free ultrafilter on $\N$. For a bounded real
sequence $(b_j)$, the number $B=\lim_{j,\mathcal U}b_j$ is
characterized by $\{j:|b_j-B|<\rho\}\in\mathcal U$ for every
$\rho>0$. Ultralimits preserve sums and inequalities. If $(y_j)$
lies in a compact metric space $K$, there is likewise a unique
$y=\lim_{j,\mathcal U}y_j\in K$, characterized by
$\{j:d(y_j,y)<\rho\}\in\mathcal U$ for every $\rho>0$.
We use this fact for bounded sets in finite-dimensional stages,
where the limits are taken in norm. Throughout this section we use
free ultrafilters on $\mathbb N$ and bounded sequences.
For a Banach space $X$, the \emph{metric type} of $(x_j)$ along
$\mathcal U$ is
\[
  \Phi_{x,\mathcal U}(z):=\lim_{j,\mathcal U}\|x_j-z\|,\qquad z\in X.
\]
When $(x_j)$ and $\mathcal U$ are clear, we simply write $\Phi$.

\begin{definition}[Radius, centre and realization of a metric type]
Let $\Phi:X\to[0,\infty)$ be a metric type. Its \emph{radius}
(or asymptotic radius for a type arising from a sequence) is
\[
 \rho(\Phi):=\inf_{z\in X}\Phi(z).
\]
An \emph{exact metric centre} of $\Phi$ is a point $v\in X$
attaining the infimum, $\Phi(v)=\rho(\Phi)$; the set of centres
may be empty or contain several points. We say that $\Phi$ is
\emph{realized} by $v\in X$ if $\Phi(z)=\|z-v\|$ for all $z\in X$.
More generally, $v$ is an \emph{exact radial centre} of $\Phi$
if there exists a number $r\ge0$ such that
\[
 \Phi(z)=\|z-v\|+r\qquad(z\in X).
\]
Then $r=\rho(\Phi)$ and $v$ is an exact metric centre.
An \emph{approximate radial centre} of $\Phi$ is a point
$v\in X$ for which there exist $0\le r\le R$ satisfying
\[
 \|z-v\|+r\le\Phi(z)\le\|z-v\|+R\qquad(z\in X).
\]
These inequalities imply $r\le\rho(\Phi)\le R$, although
$v$ need not minimize $\Phi$. The case $r=R$ gives an exact
radial centre. Both notions require the stated identities or
inequalities for every $z\in X$.
\end{definition}

\begin{proposition}[Exact radial centres in $\ell_1$-sums]
\label{prop:prel-types-l1}
Every metric type on an $\ell_1$-sum of finite-dimensional Banach
spaces has an exact radial centre.
\end{proposition}
\begin{proof}
Let $Z=(\bigoplus_{n\ge1}F_n)_{\ell_1}$, with each $F_n$
finite-dimensional, and let $\Phi_{x,\mathcal U}$ be the metric type
associated with a bounded sequence $(x_j)\subset Z$ and a free
ultrafilter $\mathcal U$. We shall find $a\in Z$ and $r\ge0$
such that
\begin{equation}\label{eq:prel-exact-type}
 \Phi_{x,\mathcal U}(z)=\|z-a\|+r\qquad(z\in Z).
\end{equation}
Let $P_N$ be the projection onto the first $N$ blocks.
For each $n$, the norm ultralimit exists:
\[
 a_n=\lim_{j,\mathcal U}x_j^{(n)}\in F_n.
\]
Put $a=(a_n)$. Since
$\sum_n\|a_n\|=\sup_N\lim_{j,\mathcal U}\|P_Nx_j\|
\le\sup_j\|x_j\|$, we have $a\in Z$; set
$r=\lim_{j,\mathcal U}\|x_j\|-\|a\|\ge0$.

First fix $z\in Z$ supported on the first $N$ blocks.
The exact $\ell_1$ decomposition gives
\begin{align*}
 \|x_j-z\|
 &=\|P_Nx_j-z\|+\|(I-P_N)x_j\|\\
 &=\|P_Nx_j-z\|+\|x_j\|-\|P_Nx_j\|.
\end{align*}
Taking ultralimits yields
\begin{align*}
 \Phi_{x,\mathcal U}(z)
 &=\|P_Na-z\|+\lim_{j,\mathcal U}\|x_j\|-\|P_Na\|\\
 &=\|a-z\|+\lim_{j,\mathcal U}\|x_j\|-\|a\|=\|a-z\|+r.
\end{align*}
The second equality uses that $z$ is supported on the first $N$
blocks and that the tails of $a$ contribute equally to $\|a-z\|$
and $\|a\|$. Finitely supported vectors are dense, and both
$\Phi_{x,\mathcal U}$ and $z\mapsto\|a-z\|+r$ are $1$-Lipschitz;
hence the identity extends to all $z\in Z$. Setting $z=a$
identifies the minimum:
\[
 r=\min_{z\in Z}\Phi_{x,\mathcal U}(z).
\]
Writing $J:Z\to Z_{\mathcal U}$ for the diagonal embedding,
we also obtain
\begin{equation}\label{eq:prel-ultrapower}
 r=\bigl\|[(x_j)]-J(a)\bigr\|_{Z_{\mathcal U}},\qquad
 \bigl\|[(x_j)]\bigr\|_{Z_{\mathcal U}}=\|a\|+r.
\end{equation}
The first identity follows from $\Phi_{x,\mathcal U}(a)=r$,
and the second from the definition of $r$.
\end{proof}

\begin{remark}[The part lying outside $Z$]
Identity \eqref{eq:prel-ultrapower} makes the geometric
interpretation precise: the ultrapower class $[(x_j)]$ can be
compared with its coordinatewise limit $a$, and $r$ measures
the distance from that class to $J(a)$. We do not assert a
canonical linear decomposition of the \emph{whole} ultrapower;
the identity refers to an individual class and its component $a$.
\end{remark}

\subsection{Property $Q$ and interlacing graphs}
For $k\ge1$, write $[\mathbb N]^k$ for the set of strictly
increasing sequences $(n_1,\ldots,n_k)$ of positive integers.
The \emph{Kalton interlacing graph} has these sequences as
vertices, with two distinct vertices adjacent when
\[
n_1\le m_1\le n_2\le m_2\le\cdots\le n_k\le m_k,
\quad\text{or the reverse interlacing inequalities.}
\]
Its shortest-path distance is denoted by $d_I$.
For an infinite $M\subset\mathbb N$, the restriction to $[M]^k$
is again a graph of the same kind after increasing enumeration.
\emph{The graph is connected}: given vertices
$A=(a_1,\ldots,a_k)$ and $B=(b_1,\ldots,b_k)$, take
$C=(c_1,\ldots,c_k)$ with $c_1>\max\{a_k,b_k\}$.
From $A$ we reach $C$ in $k$ edges via
\[
(a_2,\ldots,a_k,c_1),\ (a_3,\ldots,a_k,c_1,c_2),
\ \ldots,\ (c_1,\ldots,c_k),
\]
since consecutive vertices weakly interlace; the same argument
connects $B$ to $C$. For an infinite $M$, choose the $c_i$ in $M$.
\begin{definition}[Kalton's property $Q$]\label{def:Q}
A metric space $(Z,d)$ has property $Q$ (see also \cite{Qstab})
if there exists $C<\infty$ such that, for every $k\ge1$ and
every bounded map $f:([\mathbb N]^k,d_I)\to Z$, there is an
infinite $M\subset\mathbb N$ satisfying
\[
\sup_{\mathbf n,\mathbf m\in[M]^k}d(f(\mathbf n),f(\mathbf m))
\le C\,\Lip(f).
\]
Any such $C$ is called an \emph{admissible constant};
$Q_Z$ denotes the infimum of all admissible constants.
\end{definition}
It suffices to consider Lipschitz maps (if $\Lip(f)=\infty$,
there is nothing to prove). For fixed $k$ the graph has finite
diameter, so every Lipschitz map on it is bounded. Strictly
interlaced vertices
$n_1<m_1<n_2<m_2<\cdots<n_k<m_k$ are adjacent, and therefore
their images differ in norm by at most $\Lip(f)$.

\medskip\noindent\textbf{Trees of height $k$.}
\begin{definition}[The complete tree of height $k$]
\label{def:arbol}
For $k\ge1$ and an infinite set $L\subset\N$, put
\[
[L]^0:=\{\varnothing\},\qquad
\mathcal T_k(L):=\bigcup_{i=0}^{k}[L]^i,
\]
where $[L]^i$ consists of all strictly increasing $i$-tuples
with entries in $L$. We equip $\mathcal T_k(L)$ with the
\emph{prefix order}: $s\preceq t$ means that $s$ is the initial
segment of $t$ of length $|s|$. This is the \emph{complete rooted
tree of height $k$} used in the proofs. Its root is $\varnothing$
(level $0$), its $i$th level is $[L]^i$, and its leaves or
\emph{terminal nodes} form level $k$.

For a node $s=(s_1,\ldots,s_i)$ with $i<k$ and
$n\in L$, $n>s_i$, write
\[
s^\frown n:=(s_1,\ldots,s_i,n)\in[L]^{i+1}.
\]
For $s=\varnothing$ we set $\varnothing^\frown n=(n)$ for
any $n\in L$. These are precisely the \emph{immediate
successors} of $s$; every interior node has infinitely many
successors. A \emph{prefix} of $s$ is a node $t\preceq s$.
If $\mathbf b=(b_1,\ldots,b_k)$ is a leaf, write
$\mathbf b|i=(b_1,\ldots,b_i)$ and $\mathbf b|0=\varnothing$.
The chain
\[
\varnothing\prec\mathbf b|1\prec\cdots\prec\mathbf b|k
\]
is a \emph{complete branch}. For infinite $M\subset L$,
$\mathcal T_k(M)$ is the complete subtree with all entries
restricted to $M$.

A \emph{vector family on the tree} is a map
$s\mapsto v_s\in S_\eta$ defined on $\mathcal T_k(L)$.
In Theorem~\ref{lem:tree}, its terminal values will be
$v_{\mathbf b}=f(\mathbf b)$, and for $|s|<k$ the vectors
$v_s$ will be chosen as centres of the successor values.
For each nonroot node the \emph{increment} is
$u_{s^\frown n}:=v_{s^\frown n}-v_s$.
\end{definition}

\section{Analysis of Bourgain--Pisier spaces}
\label{sec:quotient}
\subsection{Previous results, notation and $\eta$-admissibility}
\label{subsec:lopez-admissibility}
Following the definition of López-Abad~\cite{LA}, let $T\subset B$
be normed spaces, let $E$ be another normed space, and take an
operator $u:T\to E$ with $\|u\|\le\eta<1$.
On $B\oplus_1E$ consider
\[
N_u=\{(t,-u(t)):t\in T\},\qquad
E_1=(B\oplus_1E)/N_u.
\]
The map $j:E\to E_1$ given by $j(e)=(0,e)+N_u$ is isometric.
An isometric embedding is \emph{$\eta$-admissible} if it is
isometrically equivalent to one of this form (see
\cite[Definition~2.3]{LA} and \cite[Definition~2.1]{LAarxiv}).
\paragraph{Push-out interpretation.}
The quotient above is precisely the push-out of the maps
$i:T\hookrightarrow B$ and $u:T\to E$ in the category of
normed spaces and contractive linear maps. In other words,
\[
\begin{array}{ccc}
T&\xrightarrow{u}&E\\
\big\downarrow i&&\big\downarrow j\\
B&\xrightarrow{\widetilde u}&E_1
\end{array}
\qquad
\widetilde u(b)=(b,0)+N_u
\]
is commutative, since $(t,0)-(0,u(t))=(t,-u(t))\in N_u$.
More precisely, if bounded operators $\alpha:B\to Z$ and
$\beta:E\to Z$ satisfy $\alpha i=\beta u$, there is a unique
operator $\gamma:E_1\to Z$ such that
$\gamma\widetilde u=\alpha$ and $\gamma j=\beta$;
it is given by $\gamma((b,e)+N_u)=\alpha(b)+\beta(e)$ and
satisfies $\|\gamma\|\le\max\{\|\alpha\|,\|\beta\|\}$.
This interpretation and its application to the construction
are standard. The push-outs of the auxiliary
stages should not be confused with the finite-dimensional stages
$G_s$ generating the final space.

We use the following facts from López-Abad~\cite{LA}:
\begin{enumerate}[label=(\alph*)]
\item A composition of $\eta$-admissible embeddings remains
$\eta$-admissible, \emph{with the same parameter}.
\item If $j:E\to E_1$ is $\eta$-admissible and $W\subset E$
is closed, the induced map
$\bar j:E/W\to E_1/j(W)$ is also $\eta$-admissible.
\item If $z_1+\cdots+z_p\in j(E)$, then
\begin{equation}\label{eq:source}
\sum_{i=1}^p\|z_i\|\ge
\left\|\sum_{i=1}^p z_i\right\|
+(1-\eta)\sum_{i=1}^p\dist(z_i,j(E)).
\end{equation}
\end{enumerate}
Statements (a) and (b) are contained in
\cite[Proposition~2.4(d,e)]{LA}, and (c) in
\cite[formula~(6), following Theorem~2.5]{LA}.
See also \cite[Proposition~2.3 and Lemma~1(b)]{LAarxiv} for
the corresponding admissibility arguments.

Let $(E_N)_{N\ge0}$ be an increasing chain of finite-dimensional
Banach spaces, all isometrically embedded into
\[
S_\eta=\overline{\bigcup_{N\ge0}E_N}.
\]
Assume that every inclusion $E_N\hookrightarrow E_M$, $N<M$,
is $\eta$-admissible; by (a), checking consecutive inclusions
suffices. Formula \eqref{eq:source} yields
\begin{equation}\label{eq:admiss}
\|z_1\|+\|z_2\|\ge \|z_1+z_2\|
+(1-\eta)\bigl[\dist(z_1,E_N)+\dist(z_2,E_N)\bigr],
\quad z_1+z_2\in E_N.
\tag{A}
\end{equation}
This is the admissibility inequality used in the geometric
sections. There $S_\eta$ is the \emph{abstract limit} of the
admissible chain, not the concrete Bourgain--Pisier quotient
$\mathcal L_\infty(X)/X$, which will be identified with a
subspace of $\mathcal E/X$.

\medskip\noindent\textbf{Terminology.}
A \emph{stage} means a finite-dimensional subspace $E_N$
of an increasing chain (or $F_N$ in the auxiliary quotient).
For $x\in S_\eta$ and a stage $E_N$ put
$\dist(x,E_N):=\inf_{a\in E_N}\|x-a\|$.
A \emph{minimizing sequence for $x$ relative to $E_N$} is a
sequence $(a_m)_{m\ge1}\subset E_N$ such that
$\|x-a_m\|\to\dist(x,E_N)$.
A \emph{best approximant} is a vector $e^{(N)}\in E_N$ with
$\|x-e^{(N)}\|=\dist(x,E_N)$.

\emph{Existence of a best approximant for every $x\in S_\eta$, whether or not $x$
lies in a stage.} Set $d=\dist(x,E_N)$. From the definition
of the infimum choose $a_m\in E_N$ with
$d\le\|x-a_m\|<d+1/m$. This is a bounded minimizing
sequence, since $\|a_m\|\le\|x\|+d+1$.
Because $E_N$ is finite-dimensional, some subsequence
$a_{m_j}$ converges in norm to $e^{(N)}\in E_N$.
Continuity of the norm gives
$\|x-e^{(N)}\|=\lim_j\|x-a_{m_j}\|=d$.
In particular, we do \emph{not} assume that $x\in\bigcup_NE_N$.
Best approximants need not be unique, and we may choose any of them.

\subsection{The auxiliary system}
Fix a separable infinite-dimensional $X$. In López-Abad's
system (\cite[Lemma~3.5]{LA}; see also
\cite[Lemma~3.5]{LAarxiv}) choose a countable dense set
$D=\{d_j:j\in\mathbb N\}\subset X$ and
$0<\eta<(3-\sqrt5)/4$. For every finite set
$s\subset\mathbb N$, the lemma provides spaces and maps
\[
X=\mathcal E_\varnothing,\quad
G_s\subset\mathcal E_s,\quad
j_{s,t}:\mathcal E_s\longrightarrow\mathcal E_t\quad(s\subset t),
\]
with the following properties:
\begin{enumerate}[label=(\roman*)]
\item $(\mathcal E_s,j_{s,t})$ is a coherent directed system of
isometric $\eta$-admissible embeddings;
\item $j_{s,t}(\mathcal E_s)$ has finite codimension in
$\mathcal E_t$;
\item each $G_s$ is finite-dimensional and contained in
$\mathcal E_s$, and $j_{s,t}$ maps $G_s$ into $G_t$;
\item $j_{\varnothing,s}(X_s)\subset G_s$, where
$X_s=\operatorname{span}\{d_j:j\in s\}$.
\end{enumerate}
Notice that $\mathcal E_s$ is \emph{not} finite-dimensional:
it contains an isometric copy of all of $X$.
Only its quotient by this copy of $X$ is finite-dimensional.

Choose a cofinal chain, for example $s_n=\{1,\ldots,n\}$.
Coherence allows us to regard $\mathcal E_{s_n}$ as increasing
subspaces of
\[
\mathcal E:=\overline{\bigcup_{n\ge1}\mathcal E_{s_n}},\qquad
\mathcal L_\infty(X):=\overline{\bigcup_{n\ge1}G_{s_n}}\ \subset\ \mathcal E.
\]
Property (iv) and density of $D$ imply
$X\subset\mathcal L_\infty(X)$. The larger space
$\mathcal E$ is used to study the quotient.

\medskip\noindent\textbf{Finite-dimensional quotient stages.}
For every $n$ put
\[
F_n:=\mathcal E_{s_n}/j_{\varnothing,s_n}(X),
\qquad
\bar j_{n,m}:F_n\longrightarrow F_m\quad(n<m),
\]
where $\bar j_{n,m}$ is induced by $j_{s_n,s_m}$.
By (ii) and $\mathcal E_\varnothing=X$ we have
$\dim F_n<\infty$. By \cite[Proposition~2.4(e)]{LA},
every $\bar j_{n,m}$ is an isometric $\eta$-admissible
embedding, with \emph{exactly the same} $\eta$ and no loss
as $m-n$ increases.

\smallskip
\noindent\textbf{Induced admissibility.}
Consider the representation $E_1=(B\oplus_1E)/N_u$
of an $\eta$-admissible embedding $j:E\to E_1$, and let
$W\subset E$ be closed. If $q:E\to E/W$ is the quotient,
taking the quotient first by $N_u$ and then by $j(W)$ gives the isometric identification
\[
E_1/j(W)\ \cong\
\frac{B\oplus_1(E/W)}{N_{q\circ u}},\qquad
N_{q\circ u}=\{(t,-q(u(t))):t\in T\}.
\]
Indeed, the kernel of the natural map
$B\oplus_1 E\to B\oplus_1(E/W)$ is $\{0\}\oplus W$,
and the image of $N_u$ is precisely $N_{q\circ u}$.
Since $\|q\circ u\|\le\eta$, the induced map
$E/W\to E_1/j(W)$ remains $\eta$-admissible.
Taking $E=\mathcal E_{s_n}$, $E_1=\mathcal E_{s_m}$
and $W=j_{\varnothing,s_n}(X)$ gives exactly $\bar j_{n,m}$.

\begin{proposition}[Identification of the admissible limit]
\label{prop:identify}
Let $S_\eta$ be the admissible limit of the quotient stages
$F_n=\mathcal E_{s_n}/X$. Then
\[
S_\eta\ \cong\ \mathcal E/X
\]
isometrically, and $\mathcal L_\infty(X)/X$ embeds
isometrically into $S_\eta$.
\end{proposition}
\begin{proof}
The auxiliary system comes from \cite[Lemmas~3.4--3.5]{LA},
while admissibility under quotients is covered by
\cite[Proposition~2.4(e)]{LA}. Identify the copies of $X$
in all stages using $j_{\varnothing,s_n}$.
If $e\in\mathcal E_{s_n}$, its class $e+X$ has the same norm
in $F_n$ as in $\mathcal E/X$:
\[
\|e+X\|_{F_n}
=\inf_{x\in X}\|e-x\|_{\mathcal E_{s_n}}
=\inf_{x\in X}\|e-x\|_{\mathcal E}
=\|e+X\|_{\mathcal E/X}.
\]
The middle equality uses the isometric inclusion of the stage
into $\mathcal E$, containing the whole of $X$.
Since $\bigcup_n\mathcal E_{s_n}$ is dense in $\mathcal E$,
its quotient classes are dense in $\mathcal E/X$.
Completion gives the conclusion. Finally, since
$X\subset\mathcal L_\infty(X)\subset\mathcal E$
are closed subspaces, the quotient norm of
$\mathcal L_\infty(X)/X$ is the restriction of the
quotient norm on $\mathcal E/X$.
\end{proof}

\begin{remark}[Stages of the Bourgain--Pisier space]
\label{rem:not-Y}
López-Abad's admissible embeddings are
$\mathcal E_s\hookrightarrow\mathcal E_t$; here
$\mathcal E_s$ contains $X$ and need not be finite-dimensional.
On the other hand,
$\mathcal L_\infty(X)=\overline{\bigcup_sG_s}$ is generated
by finite-dimensional stages, but the lemma only asserts
$j_{s,t}(G_s)\subset G_t$; it does \emph{not} say that
$G_s\hookrightarrow G_t$ is $\eta$-admissible.
The quotient stages $F_n=\mathcal E_{s_n}/X$, by contrast,
are finite-dimensional and $\eta$-admissible.
\end{remark}

\section{Approximate radial centres}
\label{sec:centres}
We now establish the existence of approximate radial centres.
Two consequences of the admissibility inequality (A) will be used
to control distances to finite-dimensional stages.

\subsection{Passing the admissibility inequality (A) to the limit}
\begin{lemma}\label{lem:limit}
Inequality \eqref{eq:admiss} holds for $z_1,z_2\in S_\eta$
whenever $z_1+z_2\in E_N$.
\end{lemma}
\begin{proof}
Write $e=z_1+z_2\in E_N$. Since $\bigcup_ME_M$ is dense
in $S_\eta$ and the stages are increasing, for all sufficiently
large $M$ we can choose $z_1^{(M)}\in E_M$ with
$\|z_1-z_1^{(M)}\|\to0$. Put $z_2^{(M)}=e-z_1^{(M)}$;
for $M\ge N$ this lies in $E_M$, and
$z_1^{(M)}+z_2^{(M)}=e$ \emph{exactly}.
Applying \eqref{eq:source} to $E_N\hookrightarrow E_M$ gives
\[
\|z_1^{(M)}\|+\|z_2^{(M)}\|\ge
\|e\|+(1-\eta)\sum_{i=1}^2\dist(z_i^{(M)},E_N).
\]
As $M$ increases, $z_2^{(M)}\to e-z_1=z_2$.
Both the norm and the distance function are continuous; indeed,
$|\dist(x,E_N)-\dist(y,E_N)|\le\|x-y\|$.
Passing to the limit yields \eqref{eq:admiss} on $S_\eta$.
\end{proof}

\begin{lemma}[Estimate from best approximants]\label{lem:best}
Let $x\in S_\eta$ and let $e^{(N)}\in E_N$ be a best
approximant to $x$. With
$d=\|x-e^{(N)}\|=\dist(x,E_N)$, for every $a\in E_N$,
\begin{equation}\label{eq:best}
\|a-e^{(N)}\|+(1-2\eta)d\le\|a-x\|\le\|a-e^{(N)}\|+d.
\end{equation}
\end{lemma}
\begin{proof}
The upper bound follows from the triangle inequality:
\[
\|a-x\|\le\|a-e^{(N)}\|+\|e^{(N)}-x\|
=\|a-e^{(N)}\|+d.
\]
For the lower bound we first recall the translation identity
\begin{equation}\label{eq:translation-dist}
\dist(x-e,E_N)=\dist(x,E_N)\qquad(e\in E_N).
\end{equation}
Indeed, translation by $e$ is a bijection of the linear space
$E_N$ onto itself, whence
\[
\dist(x-e,E_N)
=\inf_{a\in E_N}\|(x-e)-a\|
=\inf_{a\in E_N}\|x-(e+a)\|
=\inf_{b\in E_N}\|x-b\|.
\]
This identity requires only that $e\in E_N$, not that $e$
be a best approximant. Likewise, by symmetry and
\eqref{eq:translation-dist},
\[
\dist(a-x,E_N)=\dist(x-a,E_N)=\dist(x,E_N)
\qquad(a\in E_N).
\]
Apply \eqref{eq:admiss} to $z_1=a-x$ and
$z_2=x-e^{(N)}$. Since $z_1+z_2=a-e^{(N)}\in E_N$,
both distances in \eqref{eq:admiss} equal $d$. Hence
\[
\|a-x\|+\|x-e^{(N)}\|
\ge \|a-e^{(N)}\|+2(1-\eta)d.
\]
Using $\|x-e^{(N)}\|=d$ and
$2(1-\eta)-1=1-2\eta$ proves the claim:
\[
\|a-x\|\ge\|a-e^{(N)}\|+(1-2\eta)d.
\]
\end{proof}

\subsection{Approximate radial centres}
\begin{lemma}[Existence of approximate radial centres]\label{thm:center}
Every metric type on $S_\eta$ has an approximate radial centre.
\end{lemma}
\begin{proof}
Let $\Phi(z)=\lim_{j,\mathcal U}\|z-x_j\|$ be the type
of a bounded sequence $(x_j)\subset S_\eta$ along a free
ultrafilter $\mathcal U$. We construct $v\in S_\eta$ and
$R\ge0$ with
\begin{equation}\label{eq:center}
\|z-v\|+(1-2\eta)R\le\Phi(z)\le\|z-v\|+R
\qquad(z\in S_\eta).
\end{equation}
For each stage $E_N$, set
\[
 r_N=\lim_{j,\mathcal U}\dist(x_j,E_N).
\]
We shall also show that, if $r=\lim_{N\to\infty}r_N$, then
$(1-2\eta)R\le r\le R$.

\emph{Existence of $r$.}
Since $E_N\subset E_{N+1}$, for all $j$,
$0\le\dist(x_j,E_{N+1})\le\dist(x_j,E_N)$.
Taking ultralimits gives $0\le r_{N+1}\le r_N$.
Thus $(r_N)$ is decreasing and nonnegative, so
$r=\lim_{N\to\infty}r_N=\inf_Nr_N$ exists.

\textbf{Step 1: minimization on a fixed stage.}
Fix $N\ge0$. For each $j$ choose a best approximant
$e_j^{(N)}\in E_N$ to $x_j$. Since
$\|e_j^{(N)}\|\le 2\|x_j\|$, the sequence is bounded.
Finite-dimensionality yields a norm ultralimit
\[
c_N=\lim_{j,\mathcal U}e_j^{(N)}\in E_N.
\]
Because $e_j^{(N)}$ is a best approximant,
$\|e_j^{(N)}-x_j\|=\dist(x_j,E_N)$ for all $j$.
The reverse triangle inequality gives
\[
\left|\|c_N-x_j\|-\dist(x_j,E_N)\right|
=\left|\|c_N-x_j\|-\|e_j^{(N)}-x_j\|\right|
\le\|c_N-e_j^{(N)}\|.
\]
The last quantity has ultralimit zero. Taking scalar
ultralimits and using
$\Phi(c_N)=\lim_{j,\mathcal U}\|c_N-x_j\|$,
we obtain
\[
\boxed{\Phi(c_N)=\lim_{j,\mathcal U}\dist(x_j,E_N)=r_N.}
\]
Conversely, for $a\in E_N$ and each $j$ we have
$\|a-x_j\|\ge\dist(x_j,E_N)$, so $\Phi(a)\ge r_N$.
Therefore
\begin{equation}\label{eq:inf}
\boxed{r_N=\min_{a\in E_N}\Phi(a)=\inf_{a\in E_N}\Phi(a).}
\end{equation}

\medskip
\textbf{Step 2: a complete descent from each terminal stage $M$.}
Fix another index $M\ge0$, designating the \emph{terminal stage}
from which the descent begins. Choose $c_M\in E_M$ as in Step 1.
Apply Lemma~\ref{lem:best} to each $x_j$ and its best
approximant $e_j^{(M)}\in E_M$, then take ultralimits to get
\begin{equation}\label{eq:stageM}
\|a-c_M\|+(1-2\eta)r_M\le\Phi(a)\le\|a-c_M\|+r_M
\qquad(a\in E_M).
\end{equation}

For $0\le N\le M$, the vector $v_{N\mid M}\in E_N$ will be
obtained at stage $N$ in the descent starting from $M$. The index after
the vertical bar specifies the \emph{whole descent}, whereas
the index before it specifies the stage reached. No compatibility
among descents starting from different $M$ is assumed.
The scalar $R_{N\mid M}$ will be defined after their increments.

Put $v_{M\mid M}:=c_M$. For $i=M,M-1,\ldots,1$, choose
$v_{i-1\mid M}\in E_{i-1}$ as a best approximant to
$v_{i\mid M}$ and set
\[
u_{i\mid M}:=v_{i\mid M}-v_{i-1\mid M}.
\]
After completing the descent, define for $0\le N\le M$,
\[
\boxed{R_{N\mid M}:=r_M+\sum_{i=N+1}^{M}\|u_{i\mid M}\|\ge0.}
\]
For $N=M$ the sum is empty and $R_{M\mid M}=r_M$.
The translation identity \eqref{eq:translation-dist} applied
to $v_{i\mid M}$ and $v_{i-1\mid M}\in E_{i-1}$ yields
\[
\dist(u_{i\mid M},E_{i-1})
=\dist(v_{i\mid M},E_{i-1})=\|u_{i\mid M}\|.
\]
For $a\in E_{i-1}$ we also have
$\dist(a-v_{i\mid M},E_{i-1})=\|u_{i\mid M}\|$.
Applying \eqref{eq:admiss} to
$z_1=a-v_{i\mid M}$ and $z_2=u_{i\mid M}$, whose sum
is $a-v_{i-1\mid M}\in E_{i-1}$, and using the triangle
inequality, we obtain
\begin{equation}\label{eq:descent-one}
\|a-v_{i-1\mid M}\|+(1-2\eta)\|u_{i\mid M}\|
\le\|a-v_{i\mid M}\|
\le\|a-v_{i-1\mid M}\|+\|u_{i\mid M}\|.
\end{equation}
For $0\le N\le M$ and $a\in E_N$, inequality
\eqref{eq:descent-one} can be applied at levels
$i=N+1,\ldots,M$, since $a\in E_N\subset E_{i-1}$.
Combining these bounds and putting
$T_{N\mid M}=\sum_{i=N+1}^{M}\|u_{i\mid M}\|$ gives
\[
\|a-v_{N\mid M}\|+(1-2\eta)T_{N\mid M}
\le\|a-c_M\|
\le\|a-v_{N\mid M}\|+T_{N\mid M}.
\]
(The sum vanishes for $N=M$.) Apply \eqref{eq:stageM}
to the same $a$, which also belongs to $E_M$.
For the lower bound add $(1-2\eta)r_M$ to the leftmost
term; for the upper bound add $r_M$ to the rightmost term.
This gives
\[
\begin{aligned}
\Phi(a)&\ge\|a-c_M\|+(1-2\eta)r_M
\ge\|a-v_{N\mid M}\|+(1-2\eta)(T_{N\mid M}+r_M),\\
\Phi(a)&\le\|a-c_M\|+r_M
\le\|a-v_{N\mid M}\|+T_{N\mid M}+r_M.
\end{aligned}
\]
Since $R_{N\mid M}=r_M+T_{N\mid M}$, this establishes,
for every $a\in E_N$,
\begin{equation}\label{eq:finite}
\|a-v_{N\mid M}\|+(1-2\eta)R_{N\mid M}
\le\Phi(a)\le
\|a-v_{N\mid M}\|+R_{N\mid M}.
\end{equation}

\medskip
\textbf{Step 3: uniform bounds independent of the terminal stage.}
Fix $N\le M$. From the first inequality in
\eqref{eq:finite}, for each $a\in E_N$,
\[
\Phi(a)\ge\|a-v_{N\mid M}\|+(1-2\eta)R_{N\mid M}.
\]
Take infima over $E_N$ on both sides. Step 1 identifies the
infimum on the left as $r_N$, and $v_{N\mid M}\in E_N$
shows that the infimum of $\|a-v_{N\mid M}\|$ equals zero.
Hence
\begin{equation}\label{eq:bound-lower}
\underbrace{\inf_{a\in E_N}\Phi(a)}_{r_N}
\ge \underbrace{\inf_{a\in E_N}
\bigl(\|a-v_{N\mid M}\|+(1-2\eta)R_{N\mid M}\bigr)}_{(1-2\eta)R_{N\mid M}},
\qquad\text{that is,}\qquad (1-2\eta)R_{N\mid M}\le r_N.
\end{equation}
For the opposite bound, no infimum is needed: take
$a=v_{N\mid M}$ in the upper inequality of \eqref{eq:finite}.
The definition of $r_N$ gives
\[
r_N\le\Phi(v_{N\mid M})\le R_{N\mid M}.
\]
Finally, substitute $a=0\in E_N$ into the lower bound in
\eqref{eq:finite}:
\[
\|v_{N\mid M}\|+(1-2\eta)R_{N\mid M}\le\Phi(0).
\]
Both terms on the left are nonnegative and $1-2\eta>0$.
Therefore the following two bounds are simultaneous and
independent of both indices:
\begin{equation}\label{eq:bounds}
\boxed{\|v_{N\mid M}\|\le\Phi(0),\qquad
0\le R_{N\mid M}\le\frac{\Phi(0)}{1-2\eta}
\quad(0\le N\le M).}
\end{equation}
These bounds permit a diagonal extraction.

\medskip
\textbf{Step 4: diagonal extraction.}
We thus have the \emph{triangular family} of pairs
\[
\mathcal D_M=\bigl\{(v_{N\mid M},R_{N\mid M}):
0\le N\le M\bigr\},\qquad M=0,1,2,\ldots.
\]
We claim that there is a sequence of terminal stages
$M_j\to\infty$ such that for each $N$,
\begin{equation}\label{eq:diag-goal}
\lim_{j\to\infty}v_{N\mid M_j}=v_N\in E_N,
\qquad
\lim_{j\to\infty}R_{N\mid M_j}=R_N\in[0,\infty).
\end{equation}
For each fixed $N$, all pairs
$(v_{N\mid M},R_{N\mid M})$, $M\ge N$, lie in
\[
K_N:=\{v\in E_N:\|v\|\le\Phi(0)\}\times
\left[0,\frac{\Phi(0)}{1-2\eta}\right].
\]
This product is compact: the closed ball of the finite-dimensional
space $E_N$ and the scalar interval are compact.
\begin{enumerate}[label=\roman*)]
\item For $N=0$ the pairs $(v_{0\mid M},R_{0\mid M})$, $M\ge0$,
lie in $K_0$. Select a subsequence of terminal indices along
which both coordinates converge, and let $I_0\subset\N$ be
the infinite set of selected terminal indices.
\item Within $I_0$ choose an infinite set
$I_1\subset I_0\cap\{M\ge1\}$ on which
$(v_{1\mid M},R_{1\mid M})$ also converges.
Convergence for the pairs with first index $0$ is preserved.
Repeat for $N=2,3,\ldots$ to obtain
$I_0\supset I_1\supset I_2\supset\cdots$, so that for
each $N$, $(v_{N\mid M},R_{N\mid M})$ converges as
$M\to\infty$ along $I_N$.
\item Choose recursively $M_j\in I_j$ with
$M_1<M_2<\cdots$ and $M_j\ge j$.
For fixed $N$ and $j\ge N$ we have
$M_j\in I_j\subset I_N$; thus the stage-$N$ pairs
converge along \emph{one and the same sequence} $(M_j)$.
\end{enumerate}
This produces the limits $v_N$ and $R_N$ in
\eqref{eq:diag-goal} simultaneously for all $N$.
Passing to the limit in \eqref{eq:finite} and
\eqref{eq:bounds} with $N$ fixed gives
\begin{equation}\label{eq:finite-limit}
\|a-v_N\|+(1-2\eta)R_N\le\Phi(a)\le\|a-v_N\|+R_N
\quad(a\in E_N),\qquad (1-2\eta)R_N\le r_N\le R_N.
\end{equation}

\medskip
\textbf{Step 5: compatibility after passage to the limit.}
Compatibility was not assumed between distinct descents;
it now follows from an identity \emph{within each descent}.
Indeed, for $M>N$, the definition of $R_{N\mid M}$ gives
\[
\begin{aligned}
R_{N\mid M}-R_{N+1\mid M}
&=\left(r_M+\sum_{i=N+1}^{M}\|u_{i\mid M}\|\right)
-\left(r_M+\sum_{i=N+2}^{M}\|u_{i\mid M}\|\right)\\
&=\|u_{N+1\mid M}\|
=\|v_{N+1\mid M}-v_{N\mid M}\|.
\end{aligned}
\]
Fix $N$ and let $M=M_j\to\infty$. Since both stage-$N$
and stage-$(N+1)$ pairs converge along the sequence
constructed in Step 4, we get
\begin{equation}\label{eq:telescopic}
\boxed{R_N-R_{N+1}=\|v_{N+1}-v_N\|.}
\end{equation}
The right-hand side is nonnegative, so $(R_N)$ is decreasing
and, by \eqref{eq:finite-limit}, nonnegative.
Set $R=\lim_NR_N\ge0$. For $P>N$,
\[
\|v_P-v_N\|
\le \sum_{i=N+1}^{P}\|v_i-v_{i-1}\|
=\sum_{i=N+1}^{P}(R_{i-1}-R_i)
=R_N-R_P.
\]
Since $R_N$ converges, $(v_N)$ is Cauchy.
Let $v$ be its limit in the complete space $S_\eta$.
Letting $P\to\infty$ above gives
\begin{equation}\label{eq:tail}
\|v-v_N\|\le R_N-R.
\end{equation}

\medskip
\textbf{Step 6: conclusion.}
Fix $z\in\bigcup_NE_N$ and choose $N_0$ with $z\in E_{N_0}$.
Since the chain is increasing, the first inequality in
\eqref{eq:finite-limit}, with $a=z$, holds for every $N\ge N_0$.
Let $N\to\infty$, using $v_N\to v$ and $R_N\to R$, to get
\[
\|z-v\|+(1-2\eta)R\le\Phi(z)\le\|z-v\|+R.
\]
The type $\Phi$ is $1$-Lipschitz because all functions
$z\mapsto\|z-x_j\|$ are, and the union of the stages
is dense in $S_\eta$. By continuity, both inequalities
hold for every $z\in S_\eta$. Finally, passing to the limit
in $(1-2\eta)R_N\le r_N\le R_N$ yields
\begin{equation}\label{eq:center-radius-comparison}
(1-2\eta)R\le r\le R.
\end{equation}
\end{proof}

\section{Separation and the tree of centres}
\label{sec:tree-centres}
Approximate radial centres provide asymptotic separation of increments
from the finite-dimensional stages. This estimate allows us to
construct a tree of centres whose increments admit a simultaneous
extraction, later used to compare interlaced branches.
\subsection{Separation of successors}
The following elementary lemma uses only the triangle inequality
and the radial-centre estimate.
\begin{lemma}[Asymptotic separation]\label{lem:separation}
With the notation of Lemma~\ref{thm:center}, for every $N$,
\begin{equation}\label{eq:separation}
\lim_{j,\mathcal U}\dist(x_j-v,E_N)\ge (1-2\eta)\lim_{j,\mathcal U}\|x_j-v\|-\dist(v,E_N).
\end{equation}
\end{lemma}
\begin{proof}
We start with the triangle inequality for distances to a
subspace: for $y,z\in S_\eta$,
\[
\dist(y+z,E_N)\le \dist(y,E_N)+\dist(z,E_N).
\]
Set $y=x_j-v$ and $z=v$. Since $x_j=y+z$,
\begin{equation}\label{eq:triangle-dist}
\dist(x_j-v,E_N)
\ge\dist(x_j,E_N)-\dist(v,E_N).
\end{equation}
Taking ultralimits preserves inequalities and yields
\[
\lim_{j,\mathcal U}\dist(x_j-v,E_N)
\ge r_N-\dist(v,E_N)
\ge r-\dist(v,E_N).
\]
By \eqref{eq:center-radius-comparison}, $r\ge(1-2\eta)R$; putting
$z=v$ in its upper bound gives
\[
\lim_{j,\mathcal U}\|x_j-v\|=\Phi(v)\le R.
\]
Consequently,
\[
\lim_{j,\mathcal U}\dist(x_j-v,E_N)
\ge(1-2\eta)\lim_{j,\mathcal U}\|x_j-v\|
-\dist(v,E_N),
\]
as required.

\end{proof}

\begin{corollary}\label{cor:nearmin}
Under the assumptions and notation of Lemma~\ref{lem:separation},
if $\dist(v,E_N)<\varepsilon$, then
\begin{equation}\label{eq:nearmin}
\dist(x_j-v,E_N)\ge (1-2\eta)\|x_j-v\|-2\varepsilon
\end{equation}
for every $j$ in some set belonging to $\mathcal U$.
\end{corollary}
\begin{proof}
To prove the assertion, consider the bounded real sequence
\[
b_j=\dist(x_j-v,E_N)-(1-2\eta)\|x_j-v\|.
\]
Lemma~\ref{lem:separation} gives
$\lim_{j,\mathcal U}b_j\ge-\dist(v,E_N)>-\varepsilon$.
By the definition of an ultralimit,
$\{j:b_j>-2\varepsilon\}\in\mathcal U$;
otherwise its complement would lie in $\mathcal U$ and imply
$\lim_{j,\mathcal U}b_j\le-2\varepsilon$.
For these indices,
\[
\dist(x_j-v,E_N)\ge(1-2\eta)\|x_j-v\|-2\varepsilon.
\]
The gap between $-\varepsilon$ and $-2\varepsilon$ makes
this extraction possible without assuming the ultralimit is attained.
\end{proof}

\subsection{The tree of centres and simultaneous extraction}
Let $f:[\mathbb N]^k\to S_\eta$ be bounded.
Assign to every leaf $\mathbf b\in[\mathbb N]^k$ the value
$v_{\mathbf b}=f(\mathbf b)$, and choose the vectors at
interior nodes as approximate radial centres of the types of their
successors, by downward induction. The increments
$u_{s^\frown n}=v_{s^\frown n}-v_s$ satisfy
\[
 f(\mathbf b)=v_\varnothing+\sum_{i=1}^k u_{\mathbf b|i}.
\]
The subsequent extraction separates these increments from
finite-dimensional stages approximating previously chosen ones.

\begin{theorem}[Tree of centres]\label{lem:tree}
Let $k\ge1$, let $f:[\N]^k\to S_\eta$ be bounded, and let
$\varepsilon>0$. There are vectors
$v_s\in S_\eta$, $s\in[\N]^{\le k}$, an infinite set
$L=\{n_1<n_2<\cdots\}$ and an increasing sequence of
stage indices $p_0<p_1<\cdots$ such that, with
$u_{s^\frown n}:=v_{s^\frown n}-v_s$:
\begin{enumerate}[label=(\roman*)]
\item For each \emph{leaf}
$\mathbf b=(b_1,\ldots,b_k)\in[L]^k$,
\begin{equation}\label{eq:tree-tel}
f(\mathbf b)=v_\varnothing+\sum_{i=1}^k u_{\mathbf b|i}.
\end{equation}
\item If $s$ has entries in $\{n_1,\ldots,n_{t-1}\}$,
$|s|<k$, and $s^\frown n_t$ is a node, then
\begin{equation}\label{eq:tree-sep}
\dist(u_{s^\frown n_t},E_{p_{t-1}})
\ge (1-2\eta)\|u_{s^\frown n_t}\|-2\varepsilon.
\end{equation}
\item For every node $s^\frown n_t$ with entries in
$\{n_1,\ldots,n_t\}$,
\[
\dist(u_{s^\frown n_t},E_{p_t})<\varepsilon.
\]
\end{enumerate}
\end{theorem}
\begin{proof}
\textbf{1. Terminal values and interior centres, bottom-up.}
Fix \emph{one} free ultrafilter $\mathcal U$ on $\mathbb N$.
Every choice of an interior centre will use \emph{the same}
$\mathcal U$, rather than a separate ultrafilter for each node.
At each terminal node $s$ of length $k$, put $v_s=f(s)$.
Let
\[
B:=\sup\{\|f(\mathbf b)\|:\mathbf b\in[\N]^k\}<\infty.
\]
\emph{Base level $k$.} We have
$\|v_s\|=\|f(s)\|\le B$ for $|s|=k$.

\emph{Induction from level $i+1$ to level $i$.}
Fix $0\le i<k$ and assume all vectors at level $i+1$
have been defined with $\|v_t\|\le B$.
For each node $s$ of length $i$ put $m=\max s$,
with $\max\varnothing=0$. Rather than \emph{reindexing}
the tail $n>m$ by $j\mapsto m+j$, which might alter the
ultrafilter, use the following sequence on the \emph{original
natural indices}:
\[
x_n^{\,s}:=
\begin{cases}
v_{s^\frown n}, & n>m,\\
0, & n\le m.
\end{cases}
\]
Since $\{n>m\}$ is cofinite and belongs to $\mathcal U$,
the first $m$ terms do not affect any ultralimit.
Also $\sup_n\|x_n^{\,s}\|\le B$.
Apply Lemma~\ref{thm:center} to the metric type of this sequence, always using
the fixed ultrafilter $\mathcal U$, and choose a centre $v_s$
and a radius $R_s\ge0$. Writing
$\Phi_s(z)=\lim_{n,\mathcal U}\|z-x_n^{\,s}\|$, its lower
bound at $z=0$ gives
\[
\|v_s\|+(1-2\eta)R_s\le\Phi_s(0)
=\lim_{n,\mathcal U}\|x_n^{\,s}\|\le B.
\]
Since $1-2\eta>0$ and $R_s\ge0$, it follows that
$\|v_s\|\le B$. This construction applies to \emph{every}
node at level $i$, each depending only on the already
constructed level-$i+1$ vectors. There is no circularity and
no requirement of uniqueness.

Repeat for $i=k-1,k-2,\ldots,0$ to obtain $v_\varnothing$.
There are only $k$ levels, although each contains infinitely
many nodes. We have established
\begin{equation}\label{eq:tree-unifbound}
\sup\{\|v_s\|:s\in\mathcal T_k(\N)\}\le B,
\qquad
\|u_{s^\frown n}\|\le\|v_{s^\frown n}\|+\|v_s\|\le2B.
\end{equation}
Choosing the original indices also ensures that all sets of
good indices $A_s$ used below belong to the same ultrafilter
$\mathcal U$.

\textbf{2. Telescopic identity (i).}
Take a leaf $\mathbf b=(b_1,\ldots,b_k)$ and set
$s_0=\varnothing$, $s_i=(b_1,\ldots,b_i)$.
By definition $u_{s_i}=v_{s_i}-v_{s_{i-1}}$ for $1\le i\le k$.
Hence
\begin{align*}
 v_\varnothing+\sum_{i=1}^{k}u_{s_i}
 &=v_{s_0}+(v_{s_1}-v_{s_0})+(v_{s_2}-v_{s_1})
     +\cdots+(v_{s_k}-v_{s_{k-1}})\\
 &=v_{s_k}=f(b_1,\ldots,b_k).
\end{align*}
The intermediate terms cancel pairwise. Thus, upon subtracting
the expressions for two branches, the root $v_\varnothing$
cancels as well. This is an exact identity, with no error term.

\textbf{Plan of the extraction.}
Fix $0<\delta<\varepsilon/2$. After selecting
$n_1<\cdots<n_t$, put
\[
 H_t:=\{n_1,\ldots,n_t\}.
\]
This set contains only the \emph{integer indices} already
selected. The possible interior prefixes are the nodes of
$\mathcal T_{k-1}(H_t)$; the vectors to be approximated are
$v_s$ for $s\in\mathcal T_k(H_t)$, including terminal nodes
whenever they are present.

At step $t$ choose a stage $E_{p_t}$ such that
\[
 \dist(v_s,E_{p_t})<\delta
 \qquad(s\in\mathcal T_k(H_t)).
 \tag{$\star_t$}\label{eq:star-stage}
\]
The triangle inequality in $S_\eta/E_{p_t}$ yields, for every
node $s^\frown n_t$ with entries in $H_t$,
\begin{equation}\label{eq:two-errors}
\dist(u_{s^\frown n_t},E_{p_t})
\le\dist(v_{s^\frown n_t},E_{p_t})
 +\dist(v_s,E_{p_t})<2\delta<\varepsilon.
\end{equation}
Since $H_t$ is finite, the tree $\mathcal T_k(H_t)$ has
finitely many nodes. Thus only finitely many vectors
$\{v_s:s\in\mathcal T_k(H_t)\}$ must be approximated.
Density of $\bigcup_NE_N$ lets us approximate each one within
$\delta$ in some stage; monotonicity then supplies a single
stage large enough for all these approximations.

\textbf{3. Base case $t=1$.}
Before selecting any integer put, by convention,
\[
 H_0=\varnothing,\qquad
 \mathcal T_{k-1}(H_0)=\{\varnothing\},\qquad
 \{v_s:s\in\mathcal T_k(H_0)\}=\{v_\varnothing\}.
\]
Choose $p_0$ so that
$\dist(v_\varnothing,E_{p_0})<\delta$.
The only prefix to control at this point is $s=\varnothing$.
Apply Corollary~\ref{cor:nearmin} to $x_n=v_{(n)}$,
with centre $v_\varnothing$ and stage $E_{p_0}$.
The set
\[
 A_\varnothing:=\bigl\{n\ge1:\,
 \dist(u_{(n)},E_{p_0})\ge (1-2\eta)\|u_{(n)}\|-2\varepsilon\bigr\}
\]
belongs to $\mathcal U$ by Corollary~\ref{cor:nearmin}.
Choose $n_1\in A_\varnothing$ and put
$H_1:=H_0\cup\{n_1\}=\{n_1\}$.
Next choose $p_1>p_0$ with
$\dist(v_s,E_{p_1})<\delta$ for $s\in\mathcal T_k(H_1)$.
This establishes the required properties for $t=1$,
including \eqref{eq:two-errors}. The root is not a selected
integer and does not play the role of $n_1$.

\textbf{4. Inductive step $t\ge2$: first $n_t$, then $p_t$.}
Assume that $n_1<\cdots<n_{t-1}$ and $p_0<\cdots<p_{t-1}$
have been chosen so that $(\star_{t-1})$ holds.
For every $s\in\mathcal T_{k-1}(H_{t-1})$, apply
Corollary~\ref{cor:nearmin} to the successor sequence
$x_n=v_{s^\frown n}$, with centre $v_s$ and stage
$E_{p_{t-1}}$. By $(\star_{t-1})$,
$\dist(v_s,E_{p_{t-1}})<\delta<\varepsilon$.
Therefore
\[
 A_s:=\bigl\{n>\max s:\,
 \dist(v_{s^\frown n}-v_s,E_{p_{t-1}})
 \ge (1-2\eta)\|v_{s^\frown n}-v_s\|-2\varepsilon\bigr\}
\]
belongs to $\mathcal U$ (with the restriction $n>\max s$
vacuous for $s=\varnothing$).
Since $H_{t-1}$ is finite,
$\mathcal T_{k-1}(H_{t-1})$ contains finitely many nodes.
Moreover, \emph{all} sets $A_s$ belong to \emph{the same}
ultrafilter $\mathcal U$, fixed before any interior centres
were chosen. A finite intersection of members of a filter
remains in that filter; hence
\[
A:=\bigcap_{s\in\mathcal T_{k-1}(H_{t-1})}A_s\in\mathcal U.
\]
Let
\[
B_{t-1}:=\{n\in\mathbb N:n>n_{t-1}\}.
\]
This set is \emph{cofinite}: its complement
$\{1,\ldots,n_{t-1}\}$ is finite.
Every free ultrafilter contains all cofinite sets, since it
contains no finite set and, for each $B\subset\mathbb N$,
contains exactly one of $B$ and $\mathbb N\setminus B$.
Thus
\[
A\cap B_{t-1}\in\mathcal U.
\]
As $\varnothing\notin\mathcal U$, this intersection is
nonempty, and indeed infinite. Choose
$n_t\in A\cap B_{t-1}$. It follows that
$n_t>n_{t-1}$ and \eqref{eq:tree-sep} holds
\emph{simultaneously} for all possible earlier prefixes.

Set $H_t:=H_{t-1}\cup\{n_t\}$ and choose $p_t>p_{t-1}$
large enough for $(\star_t)$ to hold for all
$v_s$ with $s\in\mathcal T_k(H_t)$. Indeed, this tree
has finitely many nodes, so choose for each $v_s$ a stage
approximating it within $\delta$, and then choose $p_t$
larger than all those stage indices and $p_{t-1}$.
Since the stages increase and are isometrically embedded into
$S_\eta$, the approximations remain valid at the common stage.
Property (iii) follows from \eqref{eq:two-errors}; the condition
$(\star_t)$ is ready for the next step.

\textbf{5. Conclusion.}
The construction yields a strictly increasing sequence
$L=\{n_t:t\ge1\}$ and one increasing sequence of stage indices
$p_t$. Conditions (ii)--(iii) hold for \emph{every} node
formed from indices in $L$, while (i) holds before extracting $L$.
\end{proof}
\begin{remark}[Two kinds of tree vertices]
Only nodes at levels $0,\ldots,k-1$ receive vectors via the
approximate radial-centre lemma. At level $k$ the vectors are the
values of $f$. Thus ``tree of centres'' abbreviates a
construction with \emph{interior centres and terminal data}.
\end{remark}

\section{Gain and comparison of two branches}
\label{sec:two-branches}
The separation established in the previous section is now
combined with the admissibility inequality (A). We first prove
a gain estimate for partial sums, then apply it to increments
along two strictly interlaced branches. This is the decisive
estimate for property $Q$.
\begin{lemma}[Quantitative gain]\label{lem:acc}
For every stage $E\subset S_\eta$ and all $a,w\in S_\eta$,
\begin{equation}\label{eq:acc}
\|a+w\|\ge\|a\|+2(1-\eta)\dist(w,E)-\|w\|-2\dist(a,E).
\end{equation}
\end{lemma}
\begin{proof}
For $a\in E$ apply \eqref{eq:admiss} to $a+w$ and $-w$,
whose sum is $a$. Both distances to $E$ equal $\dist(w,E)$, giving
\[
\|a+w\|\ge\|a\|+2(1-\eta)\dist(w,E)-\|w\|.
\]
For arbitrary $a$ choose a best approximant $a_0\in E$.
The triangle inequality gives
\[
\|a+w\|\ge\|a_0+w\|-\|a-a_0\|,
\qquad
\|a_0\|\ge\|a\|-\|a-a_0\|.
\]
Substitution yields \eqref{eq:acc}.
A best approximant always exists for a finite-dimensional stage.
\end{proof}
We are now ready to combine the gain estimate with the
separated increments of the tree and compare two interlaced branches.
\begin{proposition}[Two strictly interlaced branches]
\label{prop:branches}
Let $\lambda_\eta=2(1-\eta)(1-2\eta)-1=1-6\eta+4\eta^2$.
Fix $k$, $f$ and $\varepsilon>0$, and let
$L=\{n_1<n_2<\cdots\}$, the stages $E_{p_t}$ and the
increments $u_s$ be those obtained from
Theorem~\ref{lem:tree} for these data.
If $\mathbf n=(n_{i_1},\ldots,n_{i_k})$ and
$\mathbf m=(n_{j_1},\ldots,n_{j_k})$ satisfy
$i_1<j_1<i_2<j_2<\cdots<i_k<j_k$, then
\begin{align}\label{eq:branch}
\|f(\mathbf n)-f(\mathbf m)\|\ge{}&\lambda_\eta\left(\sum_{r=1}^k\|u_{\mathbf n|r}\|+\sum_{r=1}^k\|u_{\mathbf m|r}\|\right)-C_{k,\eta}\varepsilon,\\
C_{k,\eta}={}&4(1-\eta)(2k-1)+(2k-1)(2k).\nonumber
\end{align}
\end{proposition}
\begin{proof}
For each $1\le r\le k$, define the signed increments $w_j$ and
associated positions $t_j$ in $L$ by
\[
\begin{aligned}
w_{2r-1}&=u_{\mathbf n|r}, & t_{2r-1}&=i_r,\\
w_{2r}&=-u_{\mathbf m|r}, & t_{2r}&=j_r.
\end{aligned}
\]
Thus $t_1<t_2<\cdots<t_{2k}$, and $n_{t_j}$ is the last
coordinate of the tree node defining the unsigned increment
underlying $w_j$. The telescopic identities in
Theorem~\ref{lem:tree}(i) give
\[
\sum_{j=1}^{2k}w_j=f(\mathbf n)-f(\mathbf m).
\]
Put $A_j=\sum_{r=1}^jw_r$. We shall apply the gain estimate
successively to $A_j=A_{j-1}+w_j$.

Fix $2\le j\le2k$ and write
\[
F_j:=E_{p_{t_j-1}}.
\]
Since $w_j$ equals $u_{s^\frown n_{t_j}}$ or its negative,
property (ii) of Theorem~\ref{lem:tree}, applied at the position
$t_j$, yields
\[
\dist(w_j,F_j)\ge(1-2\eta)\|w_j\|-2\varepsilon.
\]
Here we use $\dist(-x,F_j)=\dist(x,F_j)$ because $F_j$
is a linear subspace.

We next show that the preceding partial sum $A_{j-1}$ lies
within $(j-1)\varepsilon$ of this \emph{same} stage $F_j$.
For each $r<j$, property (iii) of
Theorem~\ref{lem:tree} gives
\[
\dist(w_r,E_{p_{t_r}})<\varepsilon.
\]
The sign causes no difficulty: an approximant to an unsigned
increment can be negated to approximate $w_r$.
Choose $b_r\in E_{p_{t_r}}$ with
$\|w_r-b_r\|<\varepsilon$.
Since $t_r<t_j$, we have $t_r\le t_j-1$, and the
stages are increasing. Consequently,
\[
E_{p_{t_r}}\subseteq E_{p_{t_j-1}}=F_j
\qquad(1\le r<j).
\]
Thus $\sum_{r=1}^{j-1}b_r\in F_j$, so
\[
\begin{aligned}
\dist(A_{j-1},F_j)
&\le\left\|A_{j-1}-\sum_{r=1}^{j-1}b_r\right\|\\
&\le\sum_{r=1}^{j-1}\|w_r-b_r\|<(j-1)\varepsilon.
\end{aligned}
\]
\emph{The decisive step: producing the gain.}
Apply Lemma~\ref{lem:acc} with $a=A_{j-1}$, $w=w_j$,
and substitute the two bounds just obtained. For
$2\le j\le2k$,
\begin{align*}
\|A_j\|
&\ge \|A_{j-1}\|+2(1-\eta)\dist(w_j,F_j)
     -\|w_j\|-2\dist(A_{j-1},F_j)\\
&\ge \|A_{j-1}\|
     +2(1-\eta)\bigl((1-2\eta)\|w_j\|-2\varepsilon\bigr)
     -\|w_j\|-2(j-1)\varepsilon\\
&=\|A_{j-1}\|+
     \underbrace{\bigl[2(1-\eta)(1-2\eta)-1\bigr]}_{\text{gain}}
     \|w_j\|
     -\varepsilon\bigl[4(1-\eta)+2(j-1)\bigr]\\
&=\|A_{j-1}\|+(1-6\eta+4\eta^2)\|w_j\|
     -\varepsilon\bigl[4(1-\eta)+2(j-1)\bigr].
\end{align*}
We now put
\[
\boxed{\lambda_\eta:=1-6\eta+4\eta^2>0.}
\]
This positivity, equivalent within the admissibility range to
$0<\eta<(3-\sqrt5)/4$, allows the norms of increments to
accumulate with a uniform gain independent of branch length.
Summing for $j=2,\ldots,2k$ and using
$\|A_1\|=\|w_1\|\ge\lambda_\eta\|w_1\|$ since
$0<\lambda_\eta\le1$, we get
\begin{align*}
\|A_{2k}\|
&\ge\lambda_\eta\sum_{j=1}^{2k}\|w_j\|
 -\varepsilon\sum_{j=2}^{2k}\bigl[4(1-\eta)+2(j-1)\bigr]\\
&=\lambda_\eta\sum_{j=1}^{2k}\|w_j\|
 -\varepsilon\bigl[4(1-\eta)(2k-1)+(2k-1)(2k)\bigr],
\end{align*}
since $2\sum_{j=2}^{2k}(j-1)=(2k-1)(2k)$.
This is \eqref{eq:branch} because
$A_{2k}=f(\mathbf n)-f(\mathbf m)$.
Notice that the restriction on $\eta$ is used
\emph{precisely} to ensure $\lambda_\eta>0$.
\end{proof}

\section{There are $\mathcal L_{\infty}$-spaces with property $Q$}
\label{sec:Q-spaces}
We first deduce property $Q$ for the admissible limit $S_\eta$
from the two-branch estimate and Ramsey's theorem. We then pass to
the Bourgain--Pisier quotient and to the associated $\mathcal L_\infty$-space.
\subsection{The space $S_{\eta}$ has property $Q$}
We shall use the following classical infinite Ramsey theorem.
\begin{theorem}[Infinite Ramsey theorem]\label{thm:ramsey}
Given $r,q\in\mathbb N$ and a colouring
$c:[\mathbb N]^r\to\{1,\ldots,q\}$, there is an infinite set
$M\subset\mathbb N$ on which $c$ is constant. The same conclusion holds
with any infinite set $L\subset\mathbb N$ in place of $\mathbb N$.
\end{theorem}
The norm of each level-$i$ increment determines a colour by partitioning
a bounded interval of $[0,\infty)$ into \emph{finitely many} intervals of
length less than $\varepsilon$. We apply Ramsey's theorem successively
for $r=1,\ldots,k$, retaining nested infinite sets
$L\supset M_1\supset\cdots\supset M_k$; once level $i$ has been
stabilised, it remains so upon passing to further subsets. This yields
\emph{one} infinite set that works simultaneously at all $k$ levels.

Applying the two-branch estimate and Ramsey's theorem to the increments
at each level will give a concentration constant independent of $k$.
\begin{theorem}[Property $Q$ of the admissible limit]\label{thm:Q}
If $0<\eta<(3-\sqrt5)/4$, then $S_\eta$ has property $Q$, with
\[Q_{S_\eta}\le\frac{1}{\lambda_\eta}
=\frac{1}{1-6\eta+4\eta^2}.\]
\end{theorem}
\begin{proof}
Fix $k\ge1$ and a Lipschitz map
$f:([\N]^k,d_I)\to S_\eta$. First suppose that $\Lip(f)=0$.
For $A,B\in[\N]^k$, choose a finite path
$A=A^{(0)},A^{(1)},\ldots,A^{(m)}=B$ in the connected interlacing
graph. Along each edge,
$\|f(A^{(j)})-f(A^{(j-1)})\|\le\Lip(f)=0$;
the triangle inequality gives $\|f(A)-f(B)\|=0$.
Thus $f$ is constant and
$\diam f([\N]^k)=0\le C\Lip(f)$ for every $C>0$, without any extraction.

Henceforth put $\ell:=\Lip(f)>0$ (reserving $L$ for the infinite set
provided by Theorem~\ref{lem:tree}), and fix $\varepsilon>0$.
Apply Theorem~\ref{lem:tree} to $f$ with precision $\varepsilon$;
write $L$, $E_{p_t}$, and $u_s$ for the resulting infinite set,
stages, and increments. We first verify that the increments are
uniformly bounded. Let
$B:=\sup_{\mathbf b\in[\N]^k}\|f(\mathbf b)\|<\infty$.
By the proof of Theorem~\ref{lem:tree}, $\|v_s\|\le B$ at every
level. Hence, for a nonroot node $s$ with parent $s^-$,
\[
\|u_s\|=\|v_s-v_{s^-}\|
\le\|v_s\|+\|v_{s^-}\|\le2B.
\]
Thus the norms of the increments belong to a bounded interval,
independently of $s$.
For each $1\le i\le k$, colour $[L]^i$ according to the interval
of a finite partition, with intervals of length less than
$\varepsilon$, containing $\|u_s\|$.
Successively applying the infinite Ramsey theorem to the $k$ levels,
and passing to nested infinite subsets, gives an infinite $M\subset L$
and numbers $a_1,\ldots,a_k\ge0$ such that
\begin{equation}\label{eq:ramsey}
|\|u_s\|-a_i|<\varepsilon\quad\text{for every }s\in[M]^i,
\ 1\le i\le k.
\end{equation}
Let $A_k=\sum_{i=1}^k a_i$. For two strictly interlaced branches
of $[M]^k$, Proposition~\ref{prop:branches} identifies the gain
coefficient
\[
\lambda_\eta=2(1-\eta)(1-2\eta)-1
=1-6\eta+4\eta^2>0.
\]
Its positivity allows us to turn the lower estimate for the
\emph{sum of norms} into an upper estimate using $\Lip(f)$.
Indeed, Proposition~\ref{prop:branches} and \eqref{eq:ramsey} give
\[\Lip(f)\ge\|f(\mathbf n)-f(\mathbf m)\|
\ge2\lambda_\eta A_k-\varepsilon\bigl[2k\lambda_\eta+C_{k,\eta}\bigr].\]
Consequently,
\[2\lambda_\eta A_k\le\Lip(f)+(2k\lambda_\eta+C_{k,\eta})\varepsilon.\]
For arbitrary branches $\mathbf p,\mathbf q\in[M]^k$, the
telescoping identity in Theorem~\ref{lem:tree} and \eqref{eq:ramsey}
yield
\[\|f(\mathbf p)-f(\mathbf q)\|
\le\sum_{i=1}^k\|u_{\mathbf p|i}\|+
\sum_{i=1}^k\|u_{\mathbf q|i}\|\le2A_k+2k\varepsilon.\]
Since $\lambda_\eta>0$, the preceding lower estimate gives
\[
2A_k\le\frac{\ell}{\lambda_\eta}
+\left(2k+\frac{C_{k,\eta}}{\lambda_\eta}\right)\varepsilon.
\]
Substitution yields, for \emph{all} $\mathbf p,\mathbf q\in[M]^k$,
\[
\|f(\mathbf p)-f(\mathbf q)\|
\le\frac{\ell}{\lambda_\eta}
+\left(4k+\frac{C_{k,\eta}}{\lambda_\eta}\right)\varepsilon.
\]
Taking the supremum, we obtain
\begin{equation}\label{eq:almostQ}
\diam f([M]^k)\le\frac{\Lip(f)}{\lambda_\eta}
+\left(4k+\frac{C_{k,\eta}}{\lambda_\eta}\right)\varepsilon.
\end{equation}
Given $C>1/\lambda_\eta$, since $\ell>0$ we choose
$\varepsilon>0$ sufficiently small at the outset (depending on
$k$, $f$, and $C$) that the second term of \eqref{eq:almostQ}
does not exceed $(C-1/\lambda_\eta)\ell$.
The case $\ell=0$ was dealt with above, so $C$ is an admissible
constant.
\end{proof}
\begin{remark}[Quantifiers]
The proof does not assert that the \emph{same} $M$ works for every
$\varepsilon>0$. This is unnecessary for property $Q$: given any
$C>1/\lambda_\eta$, choose one precision $\varepsilon$ and construct
the corresponding infinite set $M$.
\end{remark}

\begin{corollary}[Property $Q$ of the Bourgain--Pisier quotient]\label{cor:quotient}
Let $\mathcal L_\infty(X)$ be the Bourgain--Pisier space constructed
in Section~\ref{sec:quotient}, with
$0<\eta<(3-\sqrt5)/4$. Then $\mathcal L_\infty(X)/X$ has
property $Q$, and
\[
Q_{\mathcal L_\infty(X)/X}\le\frac{1}{1-6\eta+4\eta^2}.
\]
\end{corollary}
\begin{proof}
By Proposition~\ref{prop:identify}, the admissible limit
$S_\eta=\varinjlim F_n$ is isometric to $\mathcal E/X$, and
$\mathcal L_\infty(X)/X$ embeds isometrically into $S_\eta$.
The inclusions between the $F_n$ are $\eta$-admissible, so
Theorem~\ref{thm:Q} gives
\[
Q_{\mathcal E/X}\le\frac{1}{1-6\eta+4\eta^2}.
\]
Property $Q$ passes to subspaces, giving the result for
$\mathcal L_\infty(X)/X$.
\end{proof}

\subsection{Bourgain--Pisier spaces with property $Q$}
The three-space principle transfers property $Q$ from the quotient
and from $X$ to the Bourgain--Pisier space itself.
\begin{theorem}[Three-space theorem]\label{thm:three}
If $X$ is a closed subspace of a Banach space $Y$ and both $X$
and $Y/X$ have property $Q$, then $Y$ has property $Q$.
\end{theorem}
\begin{proof}
Let $Z=Y/X$ have its quotient norm, and let $q:Y\to Z$ be the
quotient map. Choose admissible constants $C_X,C_Z$ and fix $K>1$.
For every $y\in Y$ there is $x\in X$ with
$\|y-x\|\le K\|qy\|$, directly from the definition of the quotient
norm.
Fix $k$ and a Lipschitz map $f:[\N]^k\to Y$ with
$\ell=\Lip(f)$. By property $Q$ of $Z$ there is an infinite
$M_0\subset\N$ such that
$D:=\diam(qf)([M_0]^k)\le C_Z\ell$.
Fix $A_0\in[M_0]^k$ and, for every $A\in[M_0]^k$, choose
$g(A)\in X$ satisfying
\[\|f(A)-f(A_0)-g(A)\|\le KD.\]
Set $e(A)=f(A)-f(A_0)-g(A)$; thus
$\|e(A)\|\le KD$ for all $A$. If $A,B$ are adjacent
\emph{in the graph induced on $[M_0]^k$}, they are also adjacent
in the original graph, and
\[
\begin{aligned}
\|g(A)-g(B)\|
&=\|f(A)-f(B)-e(A)+e(B)\|\\
&\le\|f(A)-f(B)\|+\|e(A)\|+\|e(B)\|\\
&\le \ell+2KD.
\end{aligned}
\]
For arbitrary $A,B$, choose a shortest path of length
$m=d_I(A,B)$ in the graph $[M_0]^k$,
$A=A^{(0)},A^{(1)},\ldots,A^{(m)}=B$.
Summing over its edges gives
\[
\|g(A)-g(B)\|\le\sum_{j=1}^{m}
\|g(A^{(j)})-g(A^{(j-1)})\|
\le m(\ell+2KD)=d_I(A,B)(\ell+2KD).
\]
Hence $\Lip(g)\le \ell+2KD$ for the restricted graph metric
(which, upon renumbering, is the metric on $([\mathbb N]^k,d_I)$).
Property $Q$ of $X$ gives an infinite $M_1\subset M_0$ with
\[\diam g([M_1]^k)\le C_X(\ell+2KD).\]
For $A,B\in[M_1]^k$,
\[
 f(A)-f(B)=\bigl(g(A)-g(B)\bigr)+\bigl(e(A)-e(B)\bigr),
\]
so
\[
\|f(A)-f(B)\|\le\|g(A)-g(B)\|
+\|e(A)\|+\|e(B)\|\le\diam g([M_1]^k)+2KD.
\]
Taking suprema and using $D\le C_Z\ell$ yields
\[
\begin{aligned}
\diam f([M_1]^k)&\le C_X(\ell+2KD)+2KD\\
&\le [C_X+(2C_X+2)KC_Z]\ell.
\end{aligned}
\]
The constant is independent of $k$ and $f$.
\end{proof}

Bourgain and Delbaen constructed separable $\mathcal L_\infty$-spaces
without an isomorphic copy of $c_0$~\cite{BD1980}; the
Bourgain--Pisier construction~\cite{BP1983} gives a flexible
isometric embedding into $\mathcal L_\infty$-spaces with small
quotients. The following theorem strengthens the exclusion of $c_0$
in a nonlinear direction: even its unit ball cannot be uniformly
embedded into these examples.
\begin{theorem}[Bourgain--Pisier spaces with property $Q$]
\label{thm:bp-general}
Suppose that $X$ is a separable Banach space with property $Q$. If $\mathcal L_\infty(X)$
is a Bourgain--Pisier space constructed as in
Section~\ref{sec:quotient}, with
$0<\eta<(3-\sqrt5)/4$, then $\mathcal L_\infty(X)$ has property $Q$, with
\[
 Q_{\mathcal L_\infty(X)}
 \le Q_X+\frac{2(Q_X+1)}{\lambda_\eta},
 \qquad \lambda_\eta=1-6\eta+4\eta^2.
\]
In particular, if $X$ is reflexive, then
$Q_{\mathcal L_\infty(X)}\le 8/\lambda_\eta$.
Moreover, $B_{c_0}$ does not admit a uniform embedding into
$\mathcal L_\infty(X)$.
\end{theorem}
\begin{proof}
By Corollary~\ref{cor:quotient},
$\mathcal L_\infty(X)/X$ has property $Q$. By hypothesis, so does $X$.
Apply the three-space theorem~\ref{thm:three} to the closed
subspace $X\subset\mathcal L_\infty(X)$ and its quotient
$\mathcal L_\infty(X)/X$.
More precisely, the proof of Theorem~\ref{thm:three}, applied with
admissible constants $C_X>Q_X$ and
$C_Z>Q_{\mathcal L_\infty(X)/X}$, gives
\[
 Q_{\mathcal L_\infty(X)}\le C_X+(2C_X+2)K C_Z
 \qquad(K>1).
\]
Letting $K\downarrow1$ and the admissible constants decrease to
their infima, and using Corollary~\ref{cor:quotient}, yields the
stated bound. If $X$ is reflexive, Kalton's theorem gives $Q_X\le2$;
since $0<\lambda_\eta\le1$, this implies
$Q_{\mathcal L_\infty(X)}\le2+6/\lambda_\eta\le8/\lambda_\eta$.
The final assertion follows from Kalton's obstruction to uniform
embeddings of $B_{c_0}$ into spaces with property $Q$.
\end{proof}

\section{An application to uniform retracts in the second dual}
\label{sec:kalton-envelope}
\subsection{Property $Q$, selectors and uniform retractions}
Property $Q$ rules out even uniformly continuous maps that approximate
a section of a quotient onto $c_0$ sufficiently closely. This allows
us to apply Kalton's construction to the obstruction obtained above.
\begin{lemma}[Quantitative uniform obstruction]\label{lem:no-approx-selector}
Let $Y$ be a Banach space with property $Q$ and let
$q:Y\to c_0$ be a bounded linear operator. There is no uniformly
continuous map $s:B_{c_0}\to Y$ such that
\begin{equation}\label{eq:delta-section}
\sup_{x\in B_{c_0}}\|qs(x)-x\|_\infty<\frac12.
\end{equation}
In particular, if $q$ is surjective, it admits no uniformly
continuous section over $B_{c_0}$.
\end{lemma}
\begin{proof}
Suppose that $s$ satisfies \eqref{eq:delta-section} and choose
$0\le\delta<1/2$ with
$\sup_{x\in B_{c_0}}\|qs(x)-x\|_\infty\le\delta$.
Write
\[
\omega_s(t)=\sup\{\|s(x)-s(y)\|:\ x,y\in B_{c_0},\ \|x-y\|_\infty\le t\}.
\]
Uniform continuity gives $\omega_s(t)\to0$ as $t\downarrow0$.
Let $C_Q$ be an admissible constant for property $Q$ of $Y$.
Put $s_n=e_1+\cdots+e_n\in c_0$ and, for each $k\ge1$, define
\[
F_k(n_1,\ldots,n_k)=\frac1k\sum_{i=1}^k s_{n_i},
\qquad (n_1,\ldots,n_k)\in[\mathbb N]^k.
\]
The image lies in $B_{c_0}$. For interlaced tuples
$\mathbf n,\mathbf m$, the difference between the numbers of
indices $n_i\ge j$ and $m_i\ge j$ is at most one. Thus
\[
\|F_k(\mathbf n)-F_k(\mathbf m)\|_\infty\le\frac1k,
\qquad \Lip(s\circ F_k)\le\omega_s(1/k).
\]
Property $Q$ gives an infinite $M\subset\mathbb N$ with
\[
\operatorname{diam}(s\circ F_k)([M]^k)
\le C_Q\omega_s(1/k).
\]
Within every infinite $M$, choose
$n_1<\cdots<n_k<m_1<\cdots<m_k$.
At the coordinate $j=n_k+1$,
$F_k(\mathbf n)(j)=0$ and $F_k(\mathbf m)(j)=1$, so
\[
\operatorname{diam} F_k([M]^k)=1.
\]
The triangle inequality and boundedness of $q$ imply
\[
1\le 2\delta+\|q\| C_Q\omega_s(1/k).
\]
Letting $k\to\infty$ gives $1\le2\delta<1$, a contradiction.
\end{proof}

The following renorming construction is due to
Kalton~\cite{KaltonExamples}. For its application to
Bourgain--Pisier spaces we also use the theorem of
Argyros--Gasparis--Motakis on quotients of $\mathcal L_\infty$-spaces
\cite{AGM}.
\begin{definition}[Kalton's construction]\label{def:bpk}
Let $Y$ be a separable $\mathcal L_\infty$-space, and let
$q:Y\twoheadrightarrow c_0$ be a quotient map with
$\ker q\in\mathcal L_\infty$ and $\|q\|\le1$.
For $n\ge1$ write
\[
Y_n=(Y,\|\cdot\|_n),\qquad
\|y\|_n=\max\left\{\frac1n\|y\|_Y,\ \|qy\|_\infty\right\}.
\]
Define
\begin{equation}\label{eq:kalton-norm}
\boxed{\mathscr{K}(Y,q):=\left(\bigoplus_{n=1}^\infty Y_n\right)_{c_0}.}
\end{equation}
If the quotient map $q$ is fixed and there is no ambiguity,
we simply write $\mathscr{K}(Y)$.
\end{definition}
The space $\mathscr{K}(Y,q)$ is separable and is an
$\mathcal L_\infty$-space. Separability is immediate. For the $\mathcal L_\infty$ property we use the renorming
principle of \cite[Theorem~7.2]{KaltonExamples}:
the hypothesis $\ker q\in\mathcal L_\infty$ gives uniform local
constants for the spaces $Y_n$, and $c_0$-sums preserve the
$\mathcal L_\infty$ property.
Moreover, the first coordinate agrees isometrically with $Y$,
since $\|q\|\le1$ implies
\[
\|y\|_1=\max\{\|y\|_Y,\|qy\|_\infty\}=\|y\|_Y.
\]

Kalton's original principle is formulated in terms of the
nonexistence of a uniformly continuous selector for $q$
\cite{KaltonExamples}. Here we use the sufficient condition that
$Y$ has property $Q$, which, by Lemma~\ref{lem:no-approx-selector},
rules out even approximate selectors.
\begin{theorem}[Kalton]\label{thm:K0}
If $Y$ has property $Q$, then
$B_{\mathscr{K}(Y,q)}$ is not a uniformly continuous retract of
$B_{\mathscr{K}(Y,q)^{**}}$. In particular,
$\mathscr{K}(Y,q)$ is not a uniform retract of its bidual.
\end{theorem}
\begin{proof}
Put $K=\mathscr{K}(Y,q)$ and suppose there is a uniformly
continuous retraction
\[
\rho:B_{K^{**}}\longrightarrow B_K,
\qquad \rho(z)=z\quad(z\in B_K).
\]
For $n\ge1$, let $i_n:Y_n\to K$ and $P_n:K\to Y_n$ be
the coordinate inclusion and projection; both are contractive.
By the open mapping theorem, for each $a\in B_{c_0}$
we can choose a preimage $y_a\in Y$ of $a$, with $y_0=0$,
so that all the preimages lie in one fixed bounded ball of $Y$.

For each $n\ge1$, choose a maximal $1/n$-separated set
$D_n\subset B_{c_0}$ containing the origin. By maximality,
$D_n$ is a $1/n$-net in $B_{c_0}$.
For all sufficiently large $n$, define
\[
\varphi_n(a)=i_n(y_a)\in B_K\subset B_{K^{**}},
\qquad a\in D_n.
\]
Indeed, uniform boundedness of the preimages gives
$\|y_a\|_n=\max\{\|y_a\|_Y/n,\|a\|_\infty\}\le1$
for sufficiently large $n$.
If $a\ne b$ belong to $D_n$, then
\[
\|\varphi_n(a)-\varphi_n(b)\|_K
=\max\left\{\frac1n\|y_a-y_b\|_Y,\ \|a-b\|_\infty\right\}.
\]
As $\|a-b\|_\infty\ge1/n$ and the preimages stay in a fixed
ball, the Lipschitz constants of $\varphi_n$ are uniformly bounded.

Since $K$ is an $\mathcal L_\infty$-space, its bidual $K^{**}$
is isomorphic to an injective space by the classical
characterisation of $\mathcal L_\infty$-spaces
\cite{LindenstraussTzafririI}.
Consequently, Lipschitz maps from subsets of a metric space into
$K^{**}$ admit Lipschitz extensions with a multiplicative factor
independent of the subset. We may therefore choose extensions
\[
\widehat\varphi_n:B_{c_0}\longrightarrow K^{**},
\qquad \widehat\varphi_n|_{D_n}=\varphi_n,
\]
whose Lipschitz constants are uniformly bounded.
Choose $\ell<\infty$ with
$\Lip(\widehat\varphi_n)\le\ell$ for all the relevant indices.
The classical radial retraction
\cite{deFigueiredoKarlovitz1967,Thele1974},
\[
\mathcal R(z)=\frac{z}{\max\{1,\|z\|\}}
\quad(z\in K^{**}),
\]
takes values in $B_{K^{**}}$ and fixes that ball.
We check that it is $2$-Lipschitz. For $u,v\in K^{**}$,
put $a=\max\{1,\|u\|\}$ and $b=\max\{1,\|v\|\}$.
By symmetry suppose $a\ge b$. Then
\[
\begin{aligned}
\|\mathcal R(u)-\mathcal R(v)\|
&\le\frac{\|u-v\|}{a}
 +\frac{a-b}{ab}\,\|v\|\\
&\le\frac{\|u-v\|+(a-b)}{a}
 \le2\|u-v\|,
\end{aligned}
\]
since $\|v\|\le b$, $a\ge1$, and
$0\le a-b\le\|u-v\|$.
Thus
\[
\Psi_n=\rho\circ\mathcal R\circ\widehat\varphi_n:
B_{c_0}\longrightarrow B_K
\]
satisfies $\Psi_n(a)=\varphi_n(a)$ for $a\in D_n$.
Write the modulus of continuity of $\rho$ explicitly as
\[
\omega_\rho(t)=\sup\bigl\{\|\rho(u)-\rho(v)\|_K:
 u,v\in B_{K^{**}},\ \|u-v\|\le t\bigr\}.
\]
Uniform continuity gives $\omega_\rho(t)\to0$ as
$t\downarrow0$; furthermore,
\[
\begin{aligned}
\|\Psi_n(x)-\Psi_n(y)\|_K
&\le\omega_\rho\bigl(\|\mathcal R(\widehat\varphi_n(x))
 -\mathcal R(\widehat\varphi_n(y))\|\bigr)\\
&\le\omega_\rho(2\ell\|x-y\|_\infty).
\end{aligned}
\]
Set $s_n=P_n\Psi_n:B_{c_0}\to Y_n$. Since
$P_n$ is contractive, for $x,y\in B_{c_0}$ we have
\[
\|s_n(x)-s_n(y)\|_{Y_n}
\le\omega_\rho(2\ell\|x-y\|_\infty).
\]
In particular, all the $s_n$ have the same modulus of continuity
as maps into the spaces $Y_n$.
Since $\|z\|_Y\le n\|z\|_{Y_n}$, we also have
\[
\|s_n(x)-s_n(y)\|_Y
\le n\,\omega_\rho(2\ell\|x-y\|_\infty),
\]
so each $s_n$, regarded as a map into $Y$, is uniformly continuous.

As $qs_n(a)=a$ for $a\in D_n$, if $x\in B_{c_0}$ and
$a\in D_n$ satisfies $\|x-a\|_\infty\le1/n$, then
\[
\begin{aligned}
\|qs_n(x)-x\|_\infty
&\le\|q(s_n(x)-s_n(a))\|_\infty+\|a-x\|_\infty\\
&\le\|\Psi_n(x)-\Psi_n(a)\|_K+\frac1n\\
&\le\omega_\rho\!\left(\frac{2\ell}{n}\right)+\frac1n.
\end{aligned}
\]
Hence, for sufficiently large $n$,
\[
\sup_{x\in B_{c_0}}\|qs_n(x)-x\|_\infty
\le\omega_\rho\!\left(\frac{2\ell}{n}\right)+\frac1n
<\frac12,
\]
since $\omega_\rho(t)\to0$ as $t\downarrow0$.
This contradicts Lemma~\ref{lem:no-approx-selector}.
\end{proof}

\subsection{Applications to separable spaces}
Fix a separable infinite-dimensional reflexive Banach space $X$.
Choose a Bourgain--Pisier space
$\mathcal L_\infty(X)=\mathcal L_{\infty,\lambda}(X)$
with parameters as in Theorem~\ref{thm:bp-general},
and a quotient
$q_X:\mathcal L_\infty(X)\twoheadrightarrow c_0$
with $\mathcal L_\infty$ kernel, furnished by
Argyros--Gasparis--Motakis~\cite{AGM}.
By the convention in Definition~\ref{def:bpk}, we assume
$\|q_X\|\le1$. In this setting write
\[
\boxed{\mathscr{K}(\mathcal L_\infty(X)):=
\mathscr{K}(\mathcal L_\infty(X),q_X).}
\]
This notation suppresses both $\lambda$ and the chosen quotient;
it does not specify a canonical or unique space.

\begin{theorem}[Bourgain--Pisier--Kalton envelope]
\label{thm:reflexive-envelope}
For every separable infinite-dimensional reflexive Banach space
$X$, the space $\mathscr{K}(\mathcal L_\infty(X))$ is separable,
is an $\mathcal L_\infty$-space, contains $X$ isometrically,
and its unit ball is not a uniformly continuous retract of
its bidual unit ball.
\end{theorem}
\begin{proof}
The space $\mathscr{K}(\mathcal L_\infty(X))$ is separable
and $\mathcal L_\infty$ by the properties discussed after
Definition~\ref{def:bpk}. Kalton's theorem
\ref{thm:intro-kalton} gives property $Q$ for $X$,
and Theorem~\ref{thm:bp-general} transfers it to
$\mathcal L_\infty(X)$.
Theorem~\ref{thm:K0} rules out the uniform retraction.
The isometric embedding
$X\hookrightarrow\mathcal L_\infty(X)$, followed by the
first-coordinate inclusion into
$\mathscr{K}(\mathcal L_\infty(X),q_X)$, preserves norms.
\end{proof}

\begin{corollary}[An envelope for every separable space]
\label{cor:all-separable-envelope}
For every separable Banach space $X$ there is a separable
$\mathcal L_\infty$-space $Z$ containing $X$ linearly and
isometrically, such that $B_Z$ is not a uniformly continuous
retract of $B_{Z^{**}}$. In particular, $Z$ is not a uniform
retract of $Z^{**}$.
\end{corollary}
\begin{proof}
Let $V_X=C(B_{X^*},w^*)$. The map
$j_X:X\to V_X$ defined by $j_X(x)(x^*)=x^*(x)$ is linear and
isometric. Since $X$ is separable, $B_{X^*}$ is weak-* compact
and metrizable; hence $V_X$ is separable and $\mathcal L_\infty$.
Choose
\[
Z=V_X\oplus_\infty \mathscr{K}(\mathcal L_\infty(\ell_2)).
\]
This space is separable and $\mathcal L_\infty$, and contains
$X$ isometrically through its first summand.
Write $E=\mathscr{K}(\mathcal L_\infty(\ell_2))$.
If a uniformly continuous retraction
$\rho:B_{Z^{**}}\to B_Z$ existed, the contractive inclusion and
projection of the summand $E$ would yield
\[
B_{E^{**}}\xrightarrow{\ J_E^{**}\ }
B_{Z^{**}}\xrightarrow{\ \rho\ }
B_Z\xrightarrow{\ P_E\ }B_E.
\]
This composition would be a uniformly continuous retraction
of $B_{E^{**}}$ onto $B_E$, contradicting
Theorem~\ref{thm:reflexive-envelope} for $X=\ell_2$.
\end{proof}

\section{An abstract separation--gain criterion for property $Q$}
\label{sec:abstract-AS-AG}
We formulate an abstract criterion for property $Q$ in terms of two
geometric conditions on the same chain. In this section $X$ denotes
a separable Banach space and $(E_N)_{N\ge0}$ an increasing chain of
finite-dimensional subspaces with dense union in $X$, unless stated
otherwise. We say that a chain \emph{witnesses} $AS_\theta$ if it
satisfies the following condition; the property for a space requires
the existence of at least one witnessing chain.

\subsection{Asymptotic separation $AS_\theta$}
The threshold $\theta=1/2$ is essential for our criterion. We first
consider separation for $\theta>1/2$ and then explain why such an
improvement cannot hold universally.
\begin{definition}[Property $AS_\theta$ and witnessing chains]\label{def:AS}
Let $\tfrac12<\theta\le1$. We say that $(E_N)$
\emph{witnesses asymptotic separation $AS_\theta$} if, for every
bounded sequence $(x_j)\subset X$ and every free ultrafilter
$\mathcal U$ on $\mathbb N$, there is $v\in X$ such that
\begin{equation}\label{eq:AS-norm-bound}
\|v\|\le \sup_j\|x_j\|
\end{equation}
and, simultaneously for all $N$,
\begin{equation}\label{eq:AS-theta}
\lim_{j,\mathcal U}\dist(x_j-v,E_N)
\ge \theta\lim_{j,\mathcal U}\|x_j-v\|-\dist(v,E_N).
\end{equation}
A vector $v$ satisfying both conditions is called an
\emph{$AS_\theta$ witness for $(x_j)$, $\mathcal U$, and $(E_N)$}.
We say that $X$ satisfies $AS_\theta$ if it admits a chain
witnessing $AS_\theta$.
\end{definition}

\begin{proposition}\label{prop:AS-universal}
For every $0<\theta<1/2$, every bounded sequence $(x_j)\subset X$,
and every free ultrafilter $\mathcal U$, there is $v\in X$
satisfying \eqref{eq:AS-norm-bound} and \eqref{eq:AS-theta}
for every $N$.
\end{proposition}
\begin{proof}
Fix $0<\theta<1/2$, a bounded sequence $(x_j)\subset X$,
and a free ultrafilter $\mathcal U$. Write
\[
B:=\sup_j\|x_j\|,\qquad
\Phi(z):=\lim_{j,\mathcal U}\|x_j-z\|,\qquad
r_N:=\lim_{j,\mathcal U}\dist(x_j,E_N).
\]
The function $\Phi$ is $1$-Lipschitz. Since the stages increase,
$(r_N)$ is nonincreasing and nonnegative; set $r:=\lim_N r_N$.

For each $N$ and $j$, choose a best approximant
$e_j^{(N)}\in E_N$ to $x_j$. The sequence $(e_j^{(N)})_j$
is bounded, since $\|e_j^{(N)}\|\le2B$.
As $E_N$ is finite-dimensional,
$e^{(N)}=\lim_{j,\mathcal U}e_j^{(N)}$ exists in norm. Hence
\[
r_N=\Phi(e^{(N)})=\inf_{a\in E_N}\Phi(a).
\]
Density of $\bigcup_NE_N$ and continuity of $\Phi$ yield
\begin{equation}\label{eq:AS-universal-radius}
\boxed{r=\inf_{z\in X}\Phi(z).}
\end{equation}

\emph{Case $r>0$.} Since $2\theta<1$, there exists
$\varepsilon>0$ such that
\[
\theta(2r+3\varepsilon)<r.
\]
By \eqref{eq:AS-universal-radius}, choose $a\in X$ with
$\Phi(a)<r+\varepsilon$. By the definition of an ultralimit,
$\{j:\|x_j-a\|<\Phi(a)+\varepsilon\}$ belongs to
$\mathcal U$ and is therefore nonempty. Choose $j_0$ in this
set and put $v:=x_{j_0}$. Then $\|v\|\le B$, and
\[
\Phi(v)\le\Phi(a)+\|a-v\|
<2\Phi(a)+\varepsilon<2r+3\varepsilon.
\]
Consequently $\theta\Phi(v)<r$. For every $N,j$ the
triangle inequality for distances to subspaces gives
\[
\dist(x_j-v,E_N)\ge
\dist(x_j,E_N)-\dist(v,E_N).
\]
Taking ultralimits and recalling $r_N\ge r$, we obtain
\[
\begin{aligned}
\lim_{j,\mathcal U}\dist(x_j-v,E_N)
&\ge r_N-\dist(v,E_N)\\
&\ge r-\dist(v,E_N)\\
&\ge\theta\Phi(v)-\dist(v,E_N).
\end{aligned}
\]
Thus $v$ satisfies both separation inequalities simultaneously
at every stage.

\emph{Case $r=0$.} By \eqref{eq:AS-universal-radius}, for each
$m\ge1$ choose $a_m\in X$ with $\Phi(a_m)<1/m$.
For $m,n$,
\[
\|a_m-a_n\|
\le\lim_{j,\mathcal U}
\bigl(\|a_m-x_j\|+\|x_j-a_n\|\bigr)
=\Phi(a_m)+\Phi(a_n)<\frac1m+\frac1n.
\]
Hence $(a_m)$ is Cauchy. If $v=\lim_m a_m\in X$,
continuity of $\Phi$ gives $\Phi(v)=0$. Moreover,
\[
\|v\|\le\lim_{j,\mathcal U}\|x_j\|+\Phi(v)\le B.
\]
For every $N$,
\[
0\le\lim_{j,\mathcal U}\dist(x_j-v,E_N)
\le\lim_{j,\mathcal U}\|x_j-v\|=\Phi(v)=0.
\]
This verifies \eqref{eq:AS-theta} and completes the proof.
\end{proof}

The following example shows that strong separation alone does not
imply property $Q$.
\begin{proposition}[$c_0$ satisfies $AS_1$]
\label{prop:AS-c0}
With the coordinate stages
$E_N=\operatorname{span}\{e_1,\ldots,e_N\}$, the space $c_0$ satisfies $AS_1$.
Moreover, one may choose a witness $\widehat{x}\in c_0$ such that,
for every $N$,
\begin{equation}\label{eq:AS-c0}
\lim_{j,\mathcal U}\dist(x_j-\widehat{x},E_N)
=\lim_{j,\mathcal U}\|x_j-\widehat{x}\|_\infty.
\end{equation}
\end{proposition}
\begin{proof}
Fix a bounded sequence $(x_j)\subset c_0$ and a free
ultrafilter $\mathcal U$. Set $B=\sup_j\|x_j\|_\infty$ and
\[
x(m)=\lim_{j,\mathcal U}x_j(m),\qquad
r_N=\lim_{j,\mathcal U}\sup_{m>N}|x_j(m)|,\qquad
r=\lim_{N\to\infty}r_N.
\]
The numbers $r_N$ are nonincreasing and nonnegative, so the
limit exists. The vector $x$ belongs to $\ell_\infty$ and
$\|x\|_\infty\le B$, but need not belong to $c_0$.
For every $m>N$ we have $|x(m)|\le r_N$; hence
\begin{equation}\label{eq:c0-limsup}
\limsup_{m\to\infty}|x(m)|\le r.
\end{equation}
Define, for every $m$,
\[
\widehat{x}(m)=
\begin{cases}
0,&|x(m)|\le r,\\
\displaystyle\left(1-\frac{r}{|x(m)|}\right)x(m),&|x(m)|>r.
\end{cases}
\]
By construction,
$|x(m)-\widehat{x}(m)|\le r$, $|\widehat{x}(m)|\le|x(m)|\le B$, and
$|\widehat{x}(m)|=(|x(m)|-r)_+$.
By \eqref{eq:c0-limsup}, for every $\delta>0$, only finitely
many $m$ satisfy $|\widehat{x}(m)|>\delta$. Thus $\widehat{x}\in c_0$ and
$\|\widehat{x}\|_\infty\le B$.

Put $\delta_N=\sup_{m>N}|\widehat{x}(m)|\to0$ and
$q=\lim_{j,\mathcal U}\|x_j-\widehat{x}\|_\infty$.
Separating the first $N$ coordinates from the tail yields
\[
q\le
\max\!\left\{\max_{m\le N}|x(m)-\widehat{x}(m)|,
                r_N+\delta_N\right\}
\le\max\{r,r_N+\delta_N\}=r_N+\delta_N.
\]
Letting $N\to\infty$ gives $q\le r$.
Conversely, for every $N$,
\[
\|x_j-\widehat{x}\|_\infty\ge
\sup_{m>N}|x_j(m)|-\delta_N.
\]
Taking first the ultralimit and then the limit in $N$
gives $q\ge r$. Thus
\begin{equation}\label{eq:c0-q}
q=r.
\end{equation}
Finally, for $M>N$ and every $j$,
\[
\dist(x_j-\widehat{x},E_N)
=\sup_{m>N}|x_j(m)-\widehat{x}(m)|
\ge\sup_{m>M}|x_j(m)|-\delta_M.
\]
The ultralimit of the left-hand side is therefore at least
$r_M-\delta_M$. Letting $M\to\infty$,
\[
\lim_{j,\mathcal U}\dist(x_j-\widehat{x},E_N)\ge r.
\]
The reverse inequality follows from
$\dist(x_j-\widehat{x},E_N)\le\|x_j-\widehat{x}\|_\infty$ and
\eqref{eq:c0-q}. This proves \eqref{eq:AS-c0}, and hence $AS_1$.
\end{proof}
In particular, $AS_1$ does not imply property $Q$, since $c_0$
does not have Kalton's property $Q$~\cite{K}.

\subsection{Asymptotic gain $AG_\alpha$}
The admissibility estimate in Lemma~\ref{lem:acc} gives
\begin{equation}\label{eq:abstract-start}
 \|a+w\|\ge\|a\|+2(1-\eta)\dist(w,E)-\|w\|-2\dist(a,E).
\end{equation}
Here $\varepsilon$ controls the separation error, whereas
$\delta$ measures the distance of the preceding partial sum to the stage.
In extracting the increments we obtained
\[
 \dist(w,E)\ge\theta\|w\|-2\varepsilon,
 \qquad \dist(a,E)\le\delta.
\]
Substitution yields
\begin{align*}
\|a+w\|
&\ge\|a\|+2(1-\eta)(\theta\|w\|-2\varepsilon)
   -\|w\|-2\delta\\
&=\|a\|+\underbrace{[2(1-\eta)\theta-1]}_{\text{gain}}
 \|w\|-4(1-\eta)\varepsilon-2\delta.
\end{align*}
For a positive coefficient of $\|w\|$ we need
\begin{equation}\label{eq:abstract-condition-eta}
 2(1-\eta)\theta>1
 \quad\Longleftrightarrow\quad
 \theta>\frac{1}{2(1-\eta)}>\frac12.
\end{equation}
For $S_\eta$ we have $\theta=1-2\eta$, so
\[
 2(1-\eta)(1-2\eta)-1=1-6\eta+4\eta^2.
\]
\emph{This is precisely where we require separation to be sufficiently
stronger than $1/2$.} The correct condition is not merely
$\theta>1/2$: it depends on the strength of the gain.

This motivates the second condition in our criterion.
\begin{definition}[Asymptotic gain $AG_\alpha$]
\label{def:abstract-AG}
For $\alpha>0$, we say that the chain $(E_N)$ satisfies
\emph{asymptotic gain $AG_\alpha$} if, for every $N$,
$e\in E_N$, and $w\in X$,
\begin{equation}\label{eq:abstract-AG-restricted}
 \boxed{\|e+w\|\ge\|e\|+\alpha\dist(w,E_N)-\|w\|.}
\end{equation}
The constant $\alpha$ measures how much separation from the
stage is converted into norm growth. The definition only concerns
initial summands $e$ belonging to that stage.
\end{definition}

\noindent\textbf{Gain normalisation lemma.}
For a fixed chain $(E_N)$ and $\alpha>0$, the following are
equivalent:
\begin{enumerate}[label=(\roman*)]
\item The chain satisfies $AG_\alpha$, that is,
      \eqref{eq:abstract-AG-restricted} for $e\in E_N$.
\item For every $N$ and $a,w\in X$, the \emph{full form} holds:
\begin{equation}\label{eq:abstract-AG}
 \boxed{\|a+w\|\ge\|a\|+\alpha\dist(w,E_N)
    -\|w\|-2\dist(a,E_N).}
\end{equation}
\item There is a constant $C\ge0$ such that, for every $N$
and $a,w\in X$,
\begin{equation}\label{eq:abstract-AG-with-C}
 \|a+w\|\ge\|a\|+\alpha\dist(w,E_N)
              -\|w\|-C\dist(a,E_N).
\end{equation}
\end{enumerate}
Thus \emph{the existence of an error constant} is equivalent to
being able to take that constant equal to $2$. We do not assert
that $2$ is optimal.
\begin{proof}
The implication (ii)$\Rightarrow$(iii) is immediate. For
(iii)$\Rightarrow$(i), put $a=e\in E_N$, so the error term
vanishes. Suppose (i) and fix $a,w\in X$, $e\in E_N$.
By the triangle inequality and the restricted hypothesis,
\begin{align*}
 \|a+w\|
 &\ge\|e+w\|-\|a-e\|\\
 &\ge\|e\|+\alpha\dist(w,E_N)-\|w\|-\|a-e\|\\
 &\ge\|a\|+\alpha\dist(w,E_N)-\|w\|-2\|a-e\|.
\end{align*}
Taking the infimum over $e\in E_N$ gives
\eqref{eq:abstract-AG}, proving (i)$\Rightarrow$(ii).
\end{proof}

\noindent\textbf{Parameter restrictions.}
If $X$ is infinite-dimensional and the stages are
finite-dimensional, the restricted form with $e=0$ and
Riesz's lemma yield $\alpha\le2$.
On the other hand, in formulation
\eqref{eq:abstract-AG-with-C}, choosing $a=-w\notin E_N$
necessarily gives $C\ge\alpha$.
Neither observation makes $C$ an additional parameter of the
property: by the normalisation lemma the definition $AG_\alpha$
suffices, and in the proofs we always use the full form
\eqref{eq:abstract-AG}, with the coefficient $2$ which also
appears in the gain lemma for $S_\eta$.
For that quotient, estimate \eqref{eq:abstract-start} gives
$AG_{2(1-\eta)}$.
\par\medskip

\subsubsection{The optimal constant and an explicit family of renormings}
The following computation specifies the gain constants that can
appear in the criterion $AS_\theta+AG_\alpha$.
\begin{definition}[Optimal constant relative to a chain]
\label{def:alpha-optimal}
Let $X$ be a Banach space and $(E_N)$ an increasing dense chain
of finite-dimensional subspaces. Define
\[
 \alpha_*(X;(E_N))
 :=\sup\bigl(\{0\}\cup\{\alpha>0:
           (E_N)\text{ satisfies }AG_\alpha\}\bigr).
\]
If $X$ is infinite-dimensional, $0\le\alpha_*\le2$.
This is a quantitative constant of the norm and the chain,
not an isomorphic invariant with a fixed value.
\end{definition}

\begin{proposition}[Prescribed optimal constants]
\label{prop:ell1-renorm-alpha}
For every $0<\beta<1$ there exist a Banach space $X$
and an increasing dense chain $(E_N)$ of finite-dimensional
subspaces such that
\[
\alpha_*(X;(E_N))=\beta.
\]
\end{proposition}
\begin{proof}
For $c\ge0$, let $X_c=\ell_1$ with the norm
\[
\|x\|_c=\|x\|_1+c\|x\|_2,
\qquad E_N=\operatorname{span}\{e_1,\ldots,e_N\}.
\]
We shall prove that
\begin{equation}\label{eq:ell1-renorm-alpha}
\alpha_*(X_c;(E_N))=\frac{2}{1+c}.
\end{equation}
Given $e\in E_N$ and $w\in X_c$, write
$u=P_Nw$ and $v=(I-P_N)w$. By monotonicity of both
coordinate norms,
\[
 \dist_c(w,E_N)=\|v\|_1+c\|v\|_2.
\]
Additivity of $\ell_1$ over disjoint supports and the
$\ell_2$ triangle inequality give
\begin{align*}
 \|e+w\|_c+\|w\|_c-\|e\|_c
 &\ge 2\|v\|_1
   +c\bigl(\|e+w\|_2+\|w\|_2-\|e\|_2\bigr)\\
 &\ge 2\|v\|_1
 \ge\frac{2}{1+c}\bigl(\|v\|_1+c\|v\|_2\bigr),
\end{align*}
since $\|v\|_2\le\|v\|_1$. This proves the lower bound.

For the upper bound, fix $N\ge1$, set
$e=e_1\in E_N$ and $w_t=-e_1/2+t e_{N+1}$, $0<t<1/2$.
Then
\[
 \dist_c(w_t,E_N)=(1+c)t,
\]
whereas
\begin{align*}
 &\|e+w_t\|_c+\|w_t\|_c-\|e\|_c\\
 &\hspace{12mm}=2t+c\bigl(\sqrt{1+4t^2}-1\bigr).
\end{align*}
Validity of $AG_\alpha$ would require
\[
 \alpha(1+c)
 \le2+c\frac{\sqrt{1+4t^2}-1}{t}.
\]
Letting $t\downarrow0$ yields $\alpha\le2/(1+c)$.
Given $0<\beta<1$, choose $c=2/\beta-1>1$ and $X=X_c$.\qedhere
\end{proof}

\begin{remark}[Geometric interpretation]
The constant $\alpha_*$ may vanish (as in the obstructions we
shall study for $L_1[0,1]$ and Johnson--Lindenstrauss spaces),
or take any value in $(0,2]$ for equivalent norms on $\ell_1$.
Its positivity carries geometric information, but only
$\alpha>1$ can satisfy $\alpha\theta>1$ with $\theta\le1$,
as will be used in Theorem~\ref{thm:abstract-Q}.
The preceding examples show that $AG_\alpha$ is not an
isomorphic invariant at the level of optimal constants.
\end{remark}

\subsection{Abstract criterion: $AS_\theta$ and $AG_\alpha$}
We are now ready to state our sufficient criterion for property $Q$.
\begin{theorem}
\label{thm:abstract-Q}
Suppose that a single increasing, dense finite-dimensional
chain $(E_N)$ in $X$ satisfies the following conditions:
\begin{enumerate}[label=(\roman*)]
\item $(E_N)$ witnesses $AS_\theta$ for some
$\tfrac12<\theta\le1$;
\item $(E_N)$ satisfies $AG_\alpha$ for some $0<\alpha\le2$.
\end{enumerate}
If $\alpha\theta>1$, then $X$ has Kalton's property $Q$, with
\[
 Q_X\le\frac{1}{\alpha\theta-1}.
\]
\end{theorem}
\begin{proof}
Fix $k\ge1$ and a bounded Lipschitz map
$f:([\mathbb N]^k,d_I)\to X$; put $\ell=\Lip(f)$.
If $\ell=0$, connectivity of the graph makes $f$ constant,
and the conclusion is immediate. Henceforth suppose $\ell>0$.
Fix $\varepsilon>0$.

\paragraph{Step 1. Witnesses and tree extraction.}
Fix one free ultrafilter $\mathcal U$ for all nodes.
At the leaves $|s|=k$ put $v_s=f(s)$. At each interior node
$s$, by downward induction on $|s|$, apply $AS_\theta$ to the
successor sequence $(v_{s^\frown n})_{n>\max s}$,
retaining its natural indices (complete the initial terms
arbitrarily). This produces witnesses $v_s$.
If $B=\sup_{|s|=k}\|f(s)\|$, their norm condition gives
$\|v_s\|\le B$ at every level by downward induction.
Thus the increments satisfy
\[
 u_{s^\frown n}=v_{s^\frown n}-v_s,
 \qquad\|u_{s^\frown n}\|\le 2B.
\]
Telescoping along a complete branch
$\mathbf n=(n_1,\ldots,n_k)$ gives
\begin{equation}\label{eq:abstract-telescope}
 f(\mathbf n)=v_\varnothing+
   \sum_{i=1}^k u_{\mathbf n|i}.
\end{equation}

Extract an infinite set $L=\{n_1<n_2<\cdots\}$ and increasing
stages $E_{p_t}$ as in Theorem~\ref{lem:tree}, but use $AS_\theta$
instead of the separation inequality in Lemma~\ref{lem:separation}.
Here are the choices. Before selecting $n_t$, the stage
$E_{p_{t-1}}$ approximates to within $\varepsilon$ the
witnesses $v_s$ at the interior nodes built from previously
selected indices. For each such node, $AS_\theta$ and
$\dist(v_s,E_{p_{t-1}})<\varepsilon$ imply
\[
 \lim_{n,\mathcal U}
  \dist(v_{s^\frown n}-v_s,E_{p_{t-1}})
 \ge\theta\lim_{n,\mathcal U}\|v_{s^\frown n}-v_s\|
  -\varepsilon.
\]
By the definition of an ultralimit, the indices for which
\begin{equation}\label{eq:abstract-sep-extract}
 \dist(u_{s^\frown n},E_{p_{t-1}})
 \ge\theta\|u_{s^\frown n}\|-2\varepsilon
\end{equation}
hold on a set belonging to $\mathcal U$
(as in Corollary~\ref{cor:nearmin}).
Only finitely many prefixes have been selected; the intersection
of their corresponding sets and the cofinite tail
$\{n>n_{t-1}\}$ belongs to the same free ultrafilter and
is infinite. Choose $n_t$ from that intersection.
Next increase $p_t$ to approximate all the new witnesses and
increments to within $\varepsilon$; this can be done
simultaneously because the family is finite and the stages have
dense union. At $t=1$ begin with the root and a stage approximating
it. Thus the extraction controls the increments of all branches
ordered by their \emph{last selected index}.

\paragraph{Step 2. Gain along two interlaced branches.}
Let $L=\{n_1<n_2<\cdots\}$ be the set extracted in Step~1.
Choose two strictly interlaced branches
\[
 \mathbf n=(n_{i_1},\ldots,n_{i_k}),\qquad
 \mathbf m=(n_{j_1},\ldots,n_{j_k}),\qquad
 i_1<j_1<\cdots<i_k<j_k.
\]
For $1\le r\le k$, define their signed increments $w_j$ and
associated positions $t_j$ in $L$ by
\[
\begin{aligned}
 w_{2r-1}&=u_{\mathbf n|r}, & t_{2r-1}&=i_r,\\
 w_{2r}&=-u_{\mathbf m|r}, & t_{2r}&=j_r.
\end{aligned}
\]
Thus $t_1<\cdots<t_{2k}$, and $n_{t_j}$ is the last coordinate
of the node defining the unsigned increment underlying $w_j$.
By the telescopic identity \eqref{eq:abstract-telescope},
\[
 \sum_{j=1}^{2k}w_j=f(\mathbf n)-f(\mathbf m).
\]
Put $A_0=0$ and $A_j=\sum_{r=1}^j w_r$.
We apply the gain inequality successively to
$A_j=A_{j-1}+w_j$.

Fix $2\le j\le 2k$ and define the stage
\[
 F_j:=E_{p_{t_j-1}}.
\]
The separation estimate \eqref{eq:abstract-sep-extract}, obtained
in Step~1 at the position $t_j$, gives
\[
 \dist(w_j,F_j)\ge\theta\|w_j\|-2\varepsilon.
\]
Here the sign of $w_j$ is irrelevant, since $F_j$ is a linear
subspace. For each $r<j$, the approximation chosen in Step~1
at position $t_r$ supplies $b_r\in E_{p_{t_r}}$ such that
$\|w_r-b_r\|<\varepsilon$ (negating the approximant when needed).
Since $t_r<t_j$, the increasing stages satisfy
\[
 E_{p_{t_r}}\subseteq E_{p_{t_j-1}}=F_j
 \qquad(1\le r<j).
\]
Consequently, $\sum_{r=1}^{j-1}b_r\in F_j$, and
\[
\begin{aligned}
 \dist(A_{j-1},F_j)
 &\le\left\|A_{j-1}-\sum_{r=1}^{j-1}b_r\right\|\\
 &\le\sum_{r=1}^{j-1}\|w_r-b_r\|
 <(j-1)\varepsilon.
\end{aligned}
\]
Applying the full gain inequality \eqref{eq:abstract-AG}
with $a=A_{j-1}$, $w=w_j$ and $E=F_j$, we obtain, for
$2\le j\le2k$,
\begin{align*}
 \|A_j\|
 &\ge\|A_{j-1}\|+\alpha\dist(w_j,F_j)
       -\|w_j\|-2\dist(A_{j-1},F_j)\\
 &\ge\|A_{j-1}\|+\alpha(\theta\|w_j\|-2\varepsilon)
       -\|w_j\|-2(j-1)\varepsilon\\
 &=\|A_{j-1}\|+(\alpha\theta-1)\|w_j\|
       -[2\alpha+2(j-1)]\varepsilon.
\end{align*}
\emph{This is the decisive point:} the hypothesis
$\alpha\theta>1$ allows us to set
\[
 \boxed{\lambda:=\alpha\theta-1>0.}
\]
Each increment after the first therefore contributes at least
$\lambda\|w_j\|$, up to the indicated error.
The initial increment is treated separately:
$\|A_1\|=\|w_1\|\ge\lambda\|w_1\|$, since
$0<\lambda\le1$ follows from $\theta\le1$ and $\alpha\le2$.
Summing the inequalities for $j=2,\ldots,2k$ gives
\[
\begin{aligned}
 \|A_{2k}\|
 &\ge\lambda\sum_{j=1}^{2k}\|w_j\|
   -\varepsilon\bigl[2\alpha(2k-1)+(2k-1)(2k)\bigr]\\
 &\ge\lambda\sum_{j=1}^{2k}\|w_j\|
   -\varepsilon D(k,\alpha,2),
\end{aligned}
\]
where $D(k,\alpha,2):=4\alpha k+2k(2k-1)$ is the
same convenient upper bound used below. As $A_{2k}=f(\mathbf n)-f(\mathbf m)$,
we have proved
\begin{equation}\label{eq:abstract-two-branches}
 \boxed{
 \|f(\mathbf n)-f(\mathbf m)\|
 \ge\lambda\sum_{j=1}^{2k}\|w_j\|
      -\varepsilon D(k,\alpha,2).}
\end{equation}
The third argument of $D(k,\alpha,2)$ records the error
coefficient $2$ in the normalised full form of $AG_\alpha$;
it is not a further parameter of the property.
The error constant depends on $k$, but $k$ is fixed before
$\varepsilon$ is chosen; no height-dependent multiplicative
loss occurs.

\paragraph{Step 3. Ramsey theory and concentration.}
The increments are uniformly bounded. Thus, for each level
$i=1,\ldots,k$, we can colour $[L]^i$ by the interval of a
finite partition of $[0,2B]$ containing $\|u_s\|$.
Applying the infinite Ramsey theorem successively at every level,
we obtain an infinite $M\subset L$ and numbers
$a_1,\ldots,a_k\ge0$ such that
\[
 \left|\|u_s\|-a_i\right|<\varepsilon
 \qquad(s\in[M]^i,\ i=1,\ldots,k).
\]
Write $A=\sum_{i=1}^k a_i$. For two strictly interlaced
branches in $[M]^k$, adjacent in the interlacing graph,
\eqref{eq:abstract-two-branches} gives
\[
 \ell\ge2\lambda A-\varepsilon\,[2k\lambda+D(k,\alpha,2)],
 \qquad
 2\lambda A\le\ell+\varepsilon\,[2k\lambda+D(k,\alpha,2)].
\]
For arbitrary branches $\mathbf p,\mathbf q\in[M]^k$,
\eqref{eq:abstract-telescope} and the triangle inequality give
\[
 \|f(\mathbf p)-f(\mathbf q)\|
 \le\sum_{i=1}^k\|u_{\mathbf p|i}\|
    +\sum_{i=1}^k\|u_{\mathbf q|i}\|
 \le2A+2k\varepsilon.
\]
Solving the preceding estimate for $2A$ and substituting gives
\begin{equation}\label{eq:abstract-Qdiam}
 \boxed{\diam f([M]^k)
 \le\frac{\ell}{\lambda}+
     \varepsilon\left(4k+\frac{D(k,\alpha,2)}{\lambda}\right).}
\end{equation}
For any $D>1/\lambda$, since $\ell>0$, choose
$\varepsilon>0$ sufficiently small at the outset so that the
last term is less than $(D-1/\lambda)\ell$.
Thus $Q$ holds with every admissible constant $D>1/\lambda$.
Taking the infimum gives
$Q_X\le(\alpha\theta-1)^{-1}$.
\end{proof}

\subsubsection{Beyond the admissible setting}
The abstract theorem invites applications to spaces other than
$S_\eta$. The classical spaces $\ell_p$ are natural test cases.
They satisfy $AS_1$, but they fail the strong asymptotic gain
property. In fact, the following statements show how restrictive
asymptotic gain is: it captures a genuinely $\ell_1$-type
feature of the geometry.

Indeed, under a modest additional structural assumption,
strong asymptotic gain forces $\ell_1$ geometry.
We give two results illustrating this rigidity.
\begin{proposition}[Strong asymptotic gain forces a copy of $\ell_1$]
\label{prop:AG-contains-l1}
If an increasing dense chain $(E_N)$ of finite-dimensional
subspaces of $X$ satisfies $AG_\alpha$ for some $\alpha>1$,
then $X$ contains an isomorphic copy of $\ell_1$.
More precisely, for every $0<\delta<1-1/\alpha$ there is a
normalised sequence $(x_j)$ such that, for arbitrary scalars
$t_1,\ldots,t_m$,
\begin{equation}\label{eq:AG-l1-lower}
\left\|\sum_{j=1}^m t_jx_j\right\|
\ge\min\{1,\alpha(1-\delta)-1\}\sum_{j=1}^m|t_j|.
\end{equation}
\end{proposition}
\begin{proof}
Fix $\delta$ as in the statement and put
$\beta=\alpha(1-\delta)-1>0$.
Choose a unit vector $x_1\in E_{N_1}$ initially
($E_{N_1}\ne\{0\}$). By induction construct indices
$N_1<N_2<\cdots$ and unit vectors $x_j\in E_{N_j}$ with
\begin{equation}\label{eq:AG-riesz-stage}
  \operatorname{dist}(x_j,E_{N_{j-1}})>1-\delta
  \qquad (j\ge2).
\end{equation}
Indeed, given $N_{j-1}$, Riesz's lemma gives a unit vector
$u\in X$ at distance greater than $1-\delta/2$ from
$E_{N_{j-1}}$. By density of the union of the stages,
we can approximate $u$ and normalise in some
$E_{N_j}$, with $N_j>N_{j-1}$, preserving
\eqref{eq:AG-riesz-stage}.

Fix $m\ge2$ and arbitrary scalars $t_1,\ldots,t_m$.
For $j\ge2$ set
\[
  A_{j-1}=\sum_{i=1}^{j-1}t_ix_i,
  \qquad w_j=t_jx_j.
\]
Since the stages increase, $A_{j-1}\in E_{N_{j-1}}$ exactly;
hence $\operatorname{dist}(A_{j-1},E_{N_{j-1}})=0$.
By homogeneity of the distance,
\[
\operatorname{dist}(w_j,E_{N_{j-1}})
=|t_j|\operatorname{dist}(x_j,E_{N_{j-1}})
\ge(1-\delta)|t_j|.
\]
The full form \eqref{eq:abstract-AG} gives
\begin{align*}
 \|A_{j-1}+w_j\|
 &\ge \|A_{j-1}\|+\alpha\operatorname{dist}(w_j,E_{N_{j-1}})
       -\|w_j\|-2\operatorname{dist}(A_{j-1},E_{N_{j-1}})\\
 &\ge \|A_{j-1}\|+\bigl(\alpha(1-\delta)-1\bigr)|t_j|\\
 &=\|A_{j-1}\|+\beta|t_j|.
\end{align*}
Summing from $j=2$ to $m$ and using
$\|t_1x_1\|=|t_1|$, we obtain
\[
 \left\|\sum_{j=1}^m t_jx_j\right\|
 \ge |t_1|+\beta\sum_{j=2}^m|t_j|
 \ge \min\{1,\beta\}\sum_{j=1}^m|t_j|.
\]
The upper bound by the sum of absolute values follows from
$\|x_j\|=1$. Hence $(x_j)$ is equivalent to the canonical
basis of $\ell_1$. The error term
$2\dist(A_{j-1},E_{N_{j-1}})$ vanishes because the partial sums
were deliberately chosen inside the relevant stages.
\end{proof}

\begin{proposition}[Strong gain and a bimonotone FDD]
\label{prop:AG-bimonotone-FDD}
Let $X$ be a Banach space with a bimonotone FDD
$(F_n)_{n\ge1}$ for its original norm, and put
$E_N=F_1\oplus\cdots\oplus F_N$.
If $X$ satisfies $AG_\alpha$ with respect to this chain
for some $\alpha>1$, then the summation map induces an
isomorphism
\begin{equation}\label{eq:AG-FDD-l1}
\boxed{\quad X\simeq\left(\bigoplus_{n\ge1}F_n\right)_{\ell_1}.\quad}
\end{equation}
\end{proposition}
\begin{proof}
Let $P_N$ be the projection onto $E_N$.
Bimonotonicity gives $\|I-P_N\|\le1$.
If $w\in F_{N+1}$, then $(I-P_N)w=w$; hence for every
$e\in E_N$,
\[
 \|w\|=\|(I-P_N)(w-e)\|\le\|w-e\|.
\]
Taking $e=0$ gives the key equality
\begin{equation}\label{eq:AG-FDD-exact-dist}
\operatorname{dist}(w,E_N)=\|w\|.
\end{equation}
Given $x_i\in F_i$, for $j\ge2$ apply $AG_\alpha$ with
$a=\sum_{i<j}x_i\in E_{j-1}$ and $w=x_j$.
By \eqref{eq:AG-FDD-exact-dist},
\[
\left\|\sum_{i=1}^j x_i\right\|
\ge \left\|\sum_{i<j}x_i\right\|+(\alpha-1)\|x_j\|.
\]
Induction and the triangle inequality yield
\begin{equation}\label{eq:AG-FDD-two-sided}
\min\{1,\alpha-1\}\sum_{i=1}^m\|x_i\|
\le\left\|\sum_{i=1}^m x_i\right\|
\le\sum_{i=1}^m\|x_i\|.
\end{equation}
Since finite block sums are dense in $X$,
\eqref{eq:AG-FDD-two-sided} identifies $X$ with the
$\ell_1$-sum of its blocks.
\end{proof}

\section{An extremal example: Johnson--Lindenstrauss spaces}
\label{sec:JL-AS-AC}

\subsection{Construction of the space}
Choose a metric quotient map
\[
 T:\ell_1\longrightarrow L_1[0,1],\qquad \|T\|=1,
\]
such that for every $g\in L_1[0,1]$ there is $z\in\ell_1$
with $Tz=g$ and $\|z\|_1\le 2\|g\|_1$.
This uses the same general procedure of renorming through a
quotient map that appears in Kalton's construction, with two
differences: here we add the norm terms rather than take their
maximum, and the outer sum is $\ell_1$ rather than $c_0$.
For $0<\varepsilon<1$ put
\[
 X_\varepsilon
   =\{(\varepsilon z,Tz):z\in\ell_1\}
     \subset \ell_1\oplus_1 L_1[0,1],
 \qquad \|z\|_\varepsilon=\varepsilon\|z\|_1+\|Tz\|_1.
\]
Keep the quotient $T$ fixed, choose a sequence
$\varepsilon_n\downarrow0$, and define the Johnson--Lindenstrauss
space
\begin{equation}\label{eq:JL-space}
 \JL((\varepsilon_n))=\Bigl(\bigoplus_{n=1}^\infty X_{\varepsilon_n}\Bigr)_{\ell_1}.
\end{equation}
If the sequence is fixed, we simply write $\JL$.
The notation also depends on the chosen quotient $T$, which we
keep fixed; no canonical choice is intended.
This is the family of spaces studied by Johnson and Lindenstrauss
\cite[Proposition~1 and Theorem~1]{JL1980}:
$\JL$ is separable and $\mathcal L_1$, has the Schur and
Radon--Nikod\'ym properties, and does not embed linearly into
$\ell_1$; nor does it contain a linear copy of $L_1[0,1]$.

At the same time, each $X_{\varepsilon_n}$ is an
\emph{isometric} subspace of an $L_1$-space (because
$\ell_1\oplus_1L_1[0,1]$ is an $L_1$-space over a disjoint
union). Hence the sum \eqref{eq:JL-space} is an isometric
subspace of a suitable $L_1(\mu)$. In particular, $\JL$ has
property $Q$ by the metric stability of $L_1$~\cite{Qstab}.

The following oscillation formulas are the basic tool for the two
obstructions proved below.
\paragraph{Oscillating functions.}
Let $1/2<p<1$. Put
\[
 A_j=\{t\in[0,1]:\{2^jt\}<p\},
 \qquad r_j(t)\in\{-1,1\},
\]
where $r_j$ are the Rademacher functions and the braces in
$\{2^jt\}$ denote the fractional part.

Equivalently, $A_j$ consists of $2^j$ disjoint dyadic intervals,
\[
 A_j=\bigcup_{m=0}^{2^j-1}
 \left[\frac{m}{2^j},\frac{m+p}{2^j}\right).
\]
For example, when $p=3/4$,
$A_1=[0,3/8)\cup[1/2,7/8)$, while $A_2$ has four components
of length $3/16$. The shaded portions below show $A_1,\ldots,A_4$:

\begin{center}
\begin{tikzpicture}[x=10cm,y=0.8cm]
\foreach \j/\nn/\yy in {1/2/3.25,2/4/2.40,3/8/1.55,4/16/0.70}{
  \pgfmathtruncatemacro{\last}{\nn-1}
  \foreach \m in {0,...,\last}{
    \fill[black!42] ({\m/\nn},\yy) rectangle
                    ({(\m+0.75)/\nn},{\yy+0.40});
  }
  \draw[black!70,line width=0.35pt] (0,\yy) rectangle (1,{\yy+0.40});
  \node[anchor=east] at (-0.018,{\yy+0.20}) {$A_{\j}$};
}
\foreach \x/\lab in {0/0,0.25/{1/4},0.5/{1/2},0.75/{3/4},1/1}{
  \draw (\x,0.47)--(\x,0.58);
  \node[below] at (\x,0.42) {$\lab$};
}
\end{tikzpicture}
\end{center}
Each step doubles the number of intervals and halves their period,
while $|A_j|=p$ for every $j$. Thus $\mathbf1_{A_j}$ oscillates more
and more rapidly, even though the proportion on which it equals $1$
is unchanged.
We have weak-* convergence in $L_\infty[0,1]$:
\[
 \mathbf1_{A_j}\overset{w^*}{\longrightarrow}p\mathbf1,
 \qquad \mathbf1_{\{r_j=1\}}\overset{w^*}{\longrightarrow}
 \tfrac12\mathbf1.
\]
Indeed, it suffices to verify these limits first against step
functions constant on intervals and then approximate arbitrary
$L_1$-functions; the oscillation periods tend to zero.
The pointwise identities
\begin{align*}
 |\mathbf1_{A_j}-g|
  &=\mathbf1_{A_j}|1-g|+(1-\mathbf1_{A_j})|g|,\\
 |r_j-g|
  &=\mathbf1_{\{r_j=1\}}|1-g|
    +(1-\mathbf1_{\{r_j=1\}})|-1-g|
\end{align*}
therefore give, for every $g\in L_1[0,1]$,
\begin{align}
 \lim_j\|\mathbf1_{A_j}-g\|_1
 &=\int_0^1[p|1-g|+(1-p)|g|],
 \label{eq:JL-osc-p}\\
 \lim_j\|r_j-g\|_1
 &=\tfrac12\int_0^1(|1-g|+|-1-g|)
   =\int_0^1\max\{1,|g|\}.
 \label{eq:JL-osc-rademacher}
\end{align}
In particular, for \emph{every} finite-dimensional subspace
$F\subset L_1[0,1]$,
\begin{equation}\label{eq:JL-rademacher-distance}
 \lim_{j\to\infty}\dist(r_j,F)=1.
\end{equation}
The upper bound follows from $0\in F$ and $\|r_j\|_1=1$.
For the lower bound, choose best approximants $g_j\in F$
(with $\|g_j\|_1\le2$). Compactness of $2B_F$ gives
a convergent subsequence $g_{j_s}\to g\in F$, whence
\eqref{eq:JL-osc-rademacher} gives
\[
 \liminf_s\|r_{j_s}-g_{j_s}\|_1
 =\lim_s\|r_{j_s}-g\|_1
 =\int_0^1\max\{1,|g|\}\ge1.
\]
A subsequence of distances strictly below one would contradict
this estimate, so the distances converge to one.

\subsection{The space $\JL$ fails $AS_\theta$ for every $\theta>1/2$}
\begin{proposition}[Separation obstruction in Johnson--Lindenstrauss spaces]
\label{prop:JL-not-AS}
The space $\JL$ fails $AS_\theta$ for every
$1/2<\theta\le1$ and every increasing dense chain
of finite-dimensional subspaces.
\end{proposition}
\begin{proof}
Fix $\theta>1/2$ and choose $p\in(1/2,1)$ with
\begin{equation}\label{eq:JL-p-theta}
 2p\theta>1.
\end{equation}
For each $n$, choose liftings $z_j,z_0\in\ell_1$ with
\[
 Tz_j=\mathbf1_{A_j},\quad Tz_0=\mathbf1,
 \qquad \|z_j\|_1\le2,\quad\|z_0\|_1\le2.
\]
Let $x_j^{(n)}\in\JL$ be supported solely in block $n$,
with value $(\varepsilon_nz_j,\mathbf1_{A_j})$.
Uniformly in $j$,
\begin{equation}\label{eq:JL-sequence-B}
 B_n:=\sup_j\|x_j^{(n)}\|_\JL\le p+2\varepsilon_n.
\end{equation}
Fix a free ultrafilter $\mathcal U$ and write
\[
 \Phi_n(v)=\lim_{j,\mathcal U}\|x_j^{(n)}-v\|_\JL,
 \qquad R_n=\inf_{v\in \JL}\Phi_n(v).
\]
The vector $b_n\in\JL$ supported in block $n$ with value
$(\varepsilon_nz_0,\mathbf1)$ satisfies
\[
 \|x_j^{(n)}-b_n\|_\JL
 \le\|\mathbf1_{A_j}-\mathbf1\|_1
    +\varepsilon_n\|z_j-z_0\|_1
 \le 1-p+4\varepsilon_n.
\]
Therefore,
\begin{equation}\label{eq:JL-radius-upper}
 R_n\le 1-p+4\varepsilon_n.
\end{equation}

Now let $v\in\JL$ be a \emph{possible witness} with
$\|v\|_\JL\le B_n$. Let $g\in L_1[0,1]$ be the $L_1$
component of block $n$ of $v$. The projection from $\JL$
onto this component has norm at most one, so
$\|g\|_1\le p+2\varepsilon_n$, and
\eqref{eq:JL-osc-p} gives
\begin{align*}
 \Phi_n(v)
 &\ge \lim_j\|\mathbf1_{A_j}-g\|_1\\
 &=\int_0^1[p|1-g|+(1-p)|g|]\\
 &\ge p-(2p-1)\|g\|_1\\
 &\ge 2p(1-p)-2(2p-1)\varepsilon_n.
\end{align*}
The third line follows from $|1-t|\ge1-|t|$ and
$|t|\ge0$, using $p>1/2$.
Thus
\begin{equation}\label{eq:JL-Phi-lower}
 \Phi_n(v)\ge2p(1-p)-2(2p-1)\varepsilon_n
 \qquad(\|v\|\le B_n).
\end{equation}

Compute the minimum radius using the chain $(E_N)$.
For every $v\in\JL$ and finite-dimensional stage $E_N$,
\[
 \lim_{j,\mathcal U}\dist(x_j^{(n)}-v,E_N)
 =\inf_{e\in E_N}\Phi_n(v+e).
\]
For one inequality, choose best approximants in $E_N$ for each
$j$ and pass to their ultralimit using compactness of a bounded
ball in $E_N$; for the reverse inequality, fix $e$ first.
Since $\bigcup_NE_N$ is dense in $\JL$ and $\Phi_n$ is
continuous,
\begin{equation}\label{eq:JL-radius-stage}
 \lim_N\lim_{j,\mathcal U}\dist(x_j^{(n)}-v,E_N)
 =\inf_{e\in \JL}\Phi_n(v+e)=R_n.
\end{equation}
If $\JL$ satisfied $AS_\theta$, for each $n$ there would be
a witness $v_n$ with $\|v_n\|\le B_n$. Separation at each
stage would give
\[
 \lim_{j,\mathcal U}\dist(x_j^{(n)}-v_n,E_N)
 \ge\theta\Phi_n(v_n)-\dist(v_n,E_N).
\]
Letting $N\to\infty$ and using \eqref{eq:JL-radius-stage},
we obtain
\begin{equation}\label{eq:JL-radius-theta}
 R_n\ge\theta\Phi_n(v_n).
\end{equation}
Combining \eqref{eq:JL-radius-upper},
\eqref{eq:JL-Phi-lower}, and \eqref{eq:JL-radius-theta} yields
\[
 1-p+4\varepsilon_n
 \ge\theta\bigl[2p(1-p)-2(2p-1)\varepsilon_n\bigr].
\]
Let $n\to\infty$ and divide by $1-p>0$. The result,
$1\ge2p\theta$, contradicts \eqref{eq:JL-p-theta}.
The argument is independent of the chosen chain.\qedhere
\end{proof}

\medskip
\noindent\textbf{The case of $L_1[0,1]$.}
The same argument shows that $L_1[0,1]$ does not satisfy
$AS_\theta$ for $\theta>1/2$, with witnesses subject to
\eqref{eq:AS-norm-bound}, for any increasing dense
finite-dimensional chain. Indeed, put $x_j=\mathbf1_{A_j}$.
By \eqref{eq:JL-osc-p}, the associated type has minimum radius
$1-p$, whereas $\Phi(v)\ge2p(1-p)$ when $\|v\|_1\le p$.
For any dense chain, the limit of the distances to its stages
is the minimum radius by finite-dimensional compactness and
density, as in \eqref{eq:JL-radius-stage}.
Thus $AS_\theta$ would imply
$1-p\ge2\theta p(1-p)$, impossible when $2p\theta>1$.

\subsection{The space $\JL$ fails $AG_\alpha$ for every $\alpha>0$}
\begin{proposition}[Asymptotic gain obstruction in Johnson--Lindenstrauss spaces]
\label{prop:JL-not-AG}
The space $\JL$ fails $AG_\alpha$ for every
$\alpha>0$ and every increasing dense chain
of finite-dimensional subspaces.
\end{proposition}
\begin{proof}
Let $T:\ell_1\twoheadrightarrow L_1[0,1]$ be the quotient
used in the construction. By the open mapping theorem,
choose $K<\infty$ and $z_0,z_j\in\ell_1$ with
\[
 Tz_0=\mathbf1,\qquad Tz_j=r_j,\qquad
 \|z_0\|_1,\ \|z_j\|_1\le K\quad(j\in\N).
\]
Fix $0<t<1/2$ and, for each block $n$, define
$a_n,w_{n,j}\in\JL$ to be zero in all other coordinates,
and with components in $X_{\varepsilon_n}$ given by
\[
 a_n=(\varepsilon_n z_0,\mathbf1),\qquad
 w_{n,j}=\bigl(\varepsilon_n(-z_0/2+t z_j),
                         -\mathbf1/2+t r_j\bigr).
\]
Write the triangle defect as
$\Delta(a,w)=\|a+w\|_\JL+\|w\|_\JL-\|a\|_\JL$.
Since $r_j$ takes only the values $\pm1$ and $0<t<1/2$,
\[
 \|\mathbf1/2+t r_j\|_{L_1}
 =\|-\mathbf1/2+t r_j\|_{L_1}=1/2.
\]
Thus the $L_1$ contribution to $\Delta$ is exactly zero,
while the triangle inequality in $\ell_1$ gives
\begin{equation}\label{eq:JL-zero-defect}
 \begin{aligned}
 0\le\Delta(a_n,w_{n,j})
 &=\varepsilon_n\bigl(
   \|z_0/2+t z_j\|_1+\|-z_0/2+t z_j\|_1-\|z_0\|_1
   \bigr)\\
 &\le2Kt\varepsilon_n.
 \end{aligned}
\end{equation}
Suppose a chain $(E_N)$ satisfies $AG_\alpha$.
Fix $n$ and $\delta>0$. By density choose
$N=N(n,\delta)$ with $\dist(a_n,E_N)<\delta$.
Let $Q_n:\JL\to L_1[0,1]$ take the $L_1$ component of block
$n$; it has norm at most one. The subspace
\[
 F_N=Q_n(E_N)+\operatorname{span}\{\mathbf1\}
\]
is finite-dimensional. Since
$Q_n w_{n,j}=-\mathbf1/2+t r_j$,
\begin{align*}
 \dist(w_{n,j},E_N)
 &\ge \dist(Q_nw_{n,j},Q_nE_N)\\
 &\ge t\dist(r_j,F_N)\longrightarrow t
\end{align*}
by \eqref{eq:JL-rademacher-distance}.
Now apply the full inequality in the normalisation lemma
\eqref{eq:abstract-AG}, with its error coefficient $2$:
\[
 \alpha\dist(w_{n,j},E_N)
 \le\Delta(a_n,w_{n,j})+2\dist(a_n,E_N)
 \le 2Kt\varepsilon_n+2\delta.
\]
First let $j\to\infty$, then $\delta\downarrow0$
\emph{with $n$ fixed}; we get
$\alpha\le2K\varepsilon_n$.
Finally $n\to\infty$ and $\varepsilon_n\to0$ force
$\alpha\le0$. This order of limits requires no uniform estimate
for the stages as $n$ varies.
\end{proof}

\medskip
\noindent The sum \eqref{eq:JL-space} is an example with property
$Q$ that fails both hypotheses of Theorem~\ref{thm:abstract-Q}:
neither $AS_\theta$ with $\theta>1/2$ nor $AG_\alpha$
with $\alpha>0$ holds for any chain.
The obstruction is more striking than the one for $L_1[0,1]$,
noted after Proposition~\ref{prop:JL-not-AS}, since
$\JL$ is $\mathcal L_1$, has the Schur and Radon--Nikod\'ym
properties, and contains no copy of $L_1[0,1]$.
Nevertheless, its geometry includes, with arbitrarily small
perturbations, the phenomena responsible for the failure of
separation and gain. In the Bourgain--Pisier setting, $S_\eta$
has the Schur and Radon--Nikod\'ym properties, while the construction does not in general provide an $\mathcal L_1$ structure.

\medskip
\noindent The following table summarizes the distinctions. The $AS$ column uses the defining range $\theta>1/2$.
The positive entries record admissible parameters, without claiming optimality.
\begin{center}
\small
\renewcommand{\arraystretch}{1.25}
\begin{tabular}{|l|c|c|c|c|c|c|}
\hline
Space & $AS_\theta$ & $AG_\alpha$ & $Q$ & Schur & RNP & $\mathcal L_1$\\
\hline
$\ell_1$ & $\theta=1$ & $\alpha=2$ & Yes & Yes & Yes & Yes\\
$S_\eta$ & $\theta=1-2\eta$ & $\alpha=2(1-\eta)$ & Yes & Yes & Yes & Not guaranteed\\
$\JL$ & No & No & Yes & Yes & Yes & Yes\\
\hline
\end{tabular}
\end{center}

\section*{Author information}
Departamento de Matem\'aticas, Universidad de Extremadura,
Avenida de Elvas s/n, 06011 Badajoz, Spain.\\
\texttt{jesus@unex.es}

\section*{Declaration on the use of generative artificial intelligence}
The author used ChatGPT (OpenAI, GPT-6) during the development and
preparation of this manuscript for mathematical exploration, testing
of arguments, symbolic calculations, and assistance with exposition
and LaTeX editing. All mathematical claims, proofs, and references
were reviewed by the author, who assumes full responsibility for
the content and originality of the work.

\end{document}